\documentclass{article}
\usepackage[utf8]{inputenc}
\usepackage{amsmath,amsfonts}
\usepackage{algorithm}
\usepackage{array}
\usepackage[caption=false,font=normalsize,labelfont=sf,textfont=sf]{subfig}
\usepackage{textcomp}
\usepackage{stfloats}
\usepackage{url}
\usepackage{color}
\usepackage{verbatim}
\usepackage{graphicx}
\usepackage{cite}

\newtheorem{theorem}{Theorem}[section]
\newtheorem{lemma}[theorem]{Lemma}
\newtheorem{proposition}[theorem]{Proposition}
\newtheorem{corollary}[theorem]{Corollary}

\newtheorem{remark}[theorem]{Remark}
\newenvironment{proof}{\noindent{\em Proof:}}{\quad \hfill$\Box$\vspace{2ex}}

\def \bR {\mathbb{R}}
\def \bN {\mathbb{N}}

\def \sign  {\,{\rm sign}\,}
\def \prox  {\,{\rm Prox}\,}

\begin{document}

\title{\bf On Choosing Initial Values of Iteratively Reweighted $\ell_1$ Algorithms for the Piece-wise Exponential Penalty}
\author{Rongrong Lin, Shimin Li, and Yulan Liu\thanks{Corresponding author. School of Mathematics and Statistics, Guangdong University of Technology, Guangzhou 510520, P.R. China. Email: {\it ylliu@gdut.edu.cn}.}}

\date{ }

\maketitle



\maketitle


\begin{abstract}
Computing the proximal operator of the sparsity-promoting piece-wise exponential (PiE) penalty $1-e^{-|x|/\sigma}$ with a given shape parameter $\sigma>0$, which is treated as a popular nonconvex surrogate of $\ell_0$-norm, is fundamental in feature selection via support vector machines, image reconstruction, zero-one programming problems, compressed sensing, etc. Due to the nonconvexity of PiE, for a long time, its proximal operator is frequently evaluated via an iteratively reweighted $\ell_1$ algorithm, which substitutes PiE with its first-order approximation, however, the obtained solutions only are the critical point.  Based on the exact characterization of the proximal operator of PiE, we explore how the iteratively reweighted $\ell_1$ solution deviates from the true proximal operator in certain regions, which can be explicitly identified in terms of $\sigma$, the initial value and the regularization parameter in the definition of the proximal operator. Moreover, the initial value can be adaptively and simply chosen to ensure that the iteratively reweighted $\ell_1$ solution 
 belongs to the proximal operator of PiE.
 

\end{abstract}

\section{Introduction}

Sparse optimization problems arise in a wide range of fields, such as compressed sensing, image processing, statistics, machine learning, and among others \cite{Wang2021,Wright2022}. The so-called $\ell_0$-norm, which counts the nonzero components of a vector, is a natural penalty function to promote sparsity. Sparse solutions are more easily interpretable and generally lead to better generalization of the model performance.
Numerous studies on $\ell_0$-norm penalty optimization problem have been widely investigated in the literature \cite{Blumensath2008,Fan2020,Shen2016,Wright2022}.
However, such a nonconvex problem is NP-hard \cite{Blumensath2008}. 

To circumvent this challenge, there are a great many of $\ell_0$-norm surrogates listed in the literature \cite{LiuZhouLin,Wang2021,Zhou2023}.
The $\ell_1$-norm regularizer has received a great deal of attention for its continuity and convexity. Although it comes close to the $\ell_0$-norm, the $\ell_1$-norm frequently leads to problems with excessive punishment. To remedy this issue, nonconvex sparsity-inducing penalties have been employed to better approximate the $\ell_0$-norm and enhance sparsity, and hence have received considerable attention in sparse learning. Recent theoretical studies have shown their superiority to the convex counterparts in a variety of sparse learning settings, including the bridge $\ell_p$-norm penalty \cite{Foucart2009,LiuLin}, capped $\ell_1$ penalty \cite{Jiang2015,Zhang2010b}, transformed $\ell_1$ penalty  \cite{Zhang2017,Zhang2018}, log-sum penalty \cite{Candes2008}, minimax concave penalty \cite{Zhang2010}, smoothly clipped absolute derivation  \cite{Fan2001}, the difference of $\ell_1$- and $\ell_2$-norms  \cite{Lou2018,Yin2015}, the ratio of $\ell_1$- and $\ell_2$-norms  \cite{Tao2022,Yin2014}, Weibull penalty \cite{Zhou2022}, generalized error functions \cite{Guo2021,Zhou2023}, $p$-th power of the $\ell_1$-norm \cite{Prater2023},  piece-wise exponential function (PiE) in \cite{Bradley1998,Le2015,Mangasarian1996}, and among others. 
To address the nonconvex and possibly nonsmooth problems, a proximal algorithm is commonly used \cite{Beck2017}. The proximal operator \cite{Beck2017} of a function $\varphi:\bR\to\bR$ at $\tau\in\bR$ with the regularization parameter $\lambda>0$ is defined by 
$$
\prox_{\lambda \varphi}(\tau):=\arg\min_{x\in\bR} \Big\{\lambda \varphi(x)+\frac{1}{2}(x-\tau)^2 \Big\}.
$$
Characterizing the proximal operator of a function is crucial  to the proximal algorithm. 
However, such a proximal operator does not always have a closed form or
is computationally challenging to solve due to the nonconvex and nonsmooth nature of the sparsity-inducing penalty. A popular method for handling this issue is the iteratively reweighted algorithm, which approximates the nonconvex and nonsmooth problem by a sequence of trackable convex subproblems.
Zou and Li \cite{Zou2008} devised a local linear approximation, which can be treated as a special case of the iteratively reweighted $\ell_1$ (IRL1) minimization method proposed by Cand\'{e}s,
Wakin, and Boyd \cite{Candes2008}. 
The IRL1 algorithm can be unified under a majorization-minimization framework \cite{Ochs2015}.
Later, the IRL1 algorithm for optimization problems with general nonconvex and nonsmooth sparsity-inducing terms was explored in \cite{Wang2021},
and its  global and local convergence analysis for the $\ell_p$-norm regularized model were studied in \cite{WangZeng2021} and \cite{Wang2023}, respectively.

In this paper, we focus on the PiE function.
The PiE function $f_{\sigma}:\mathbb{R}\to \mathbb{R}$ with a shape parameter $\sigma>0$, defined by
\begin{equation}\label{PiE}
f_{\sigma}(x)=1-e^{-\frac{|x|}{\sigma}},  \mbox{ for any } x\in \mathbb{R},
\end{equation}
is one of the nonconvex surrogates of the $\ell_0$-norm.
It is also called an exponential-type penalty \cite{Gao2011,Le2015,Wang2021,Zhou2022} or a Laplacian function \cite[(16)]{Trzasko2008}, which has been successfully applied in the support vector machines \cite{Bradley1998,Fung2002}, zero-one programming problems \cite{Lucidi2010, Rinaldi2009}, image reconstruction \cite{Trzasko2008,Zhang2018}, compressed sensing \cite{Chen2014,Le2015,Malek2016}, and the low-rank matrix completion \cite{Yan2022}, etc. 
Due to the nonconvexity of PiE, 
for a very long time, the IRL1 algorithm was adopted in a large volume of references to approximate the proximal operator of PiE \cite{Bradley1998,Yan2022,Zhou2022,Zhou2023}. 
Recently, the IRL1 algorithm for computing the proximal operator of PiE was adopted in \cite[(3.19)]{Yan2022} for matrix completion. However,
the expression of the proximal operator $\prox_{\lambda f_{\sigma}}$ for PiE  was originally and partially studied by Malek-Mohammadi et al\cite{Malek2016} in 2016 and then systematically explored by Liu, Zhou, and Lin \cite{LiuZhouLin} using the Lambert W function. 
Motivated by the analysis between the IRL1 algorithm solution for the log-sum penalty and its proximal operator in \cite{Prater2022}, we will explore
 the relation between the IRL1 algorithm solution and the proximal operator for PiE and then provide how to select a suitable initial point in the  IRL1 algorithm to ensure that the IRL1 solution is consistent with the proximal operator of PiE.

The remainder of the paper is outlined as follows: In Section \ref{Section2}, we recall the existing characterizations for $\prox_{\lambda f_{\sigma}}$ by utilizing the Lambert W function. With this, we show in Theorems \ref{MainTh1} and \ref{MainTh2} of Section \ref{Section3} that the iteratively reweighted $\ell_1$ solution does not belong to the  proximal operator of PiE in certain regions, which can be explicitly determined in terms of $\sigma$, the initial value, and the regularization parameter $\lambda$, as shown in Fig. \ref{IRL1} later. To remedy this issue, the initial value is set adaptively, as in Theorems \ref{Theorem3} and \ref{Theorem4}, to ensure that the IRL1 solution belongs to the proximal operator of PiE.
Some necessary lemmas and the proofs of Theorems \ref{MainTh1} and \ref{MainTh2} are presented in Section \ref{Section4}.
 Some conclusions are made in the final section.

\section{Existing characterizations for $\prox_{\lambda f_{\sigma}}$}\label{Section2}

\begin{figure}[!tbp]
\begin{center}
\includegraphics[scale=0.4]{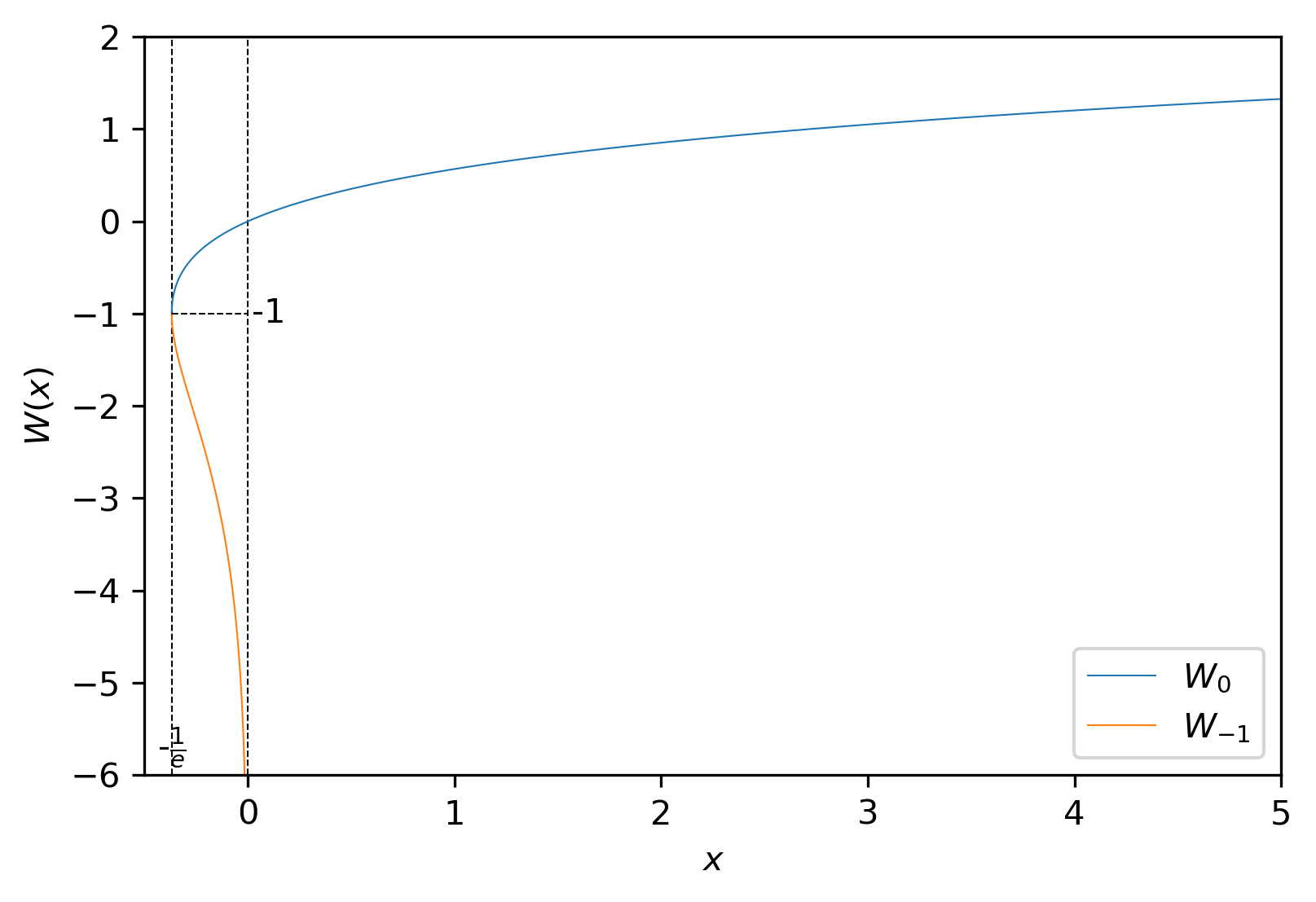}
\caption{Two main branches of the Lambert W function.}
\label{LambertW}
\end{center}
\end{figure}

Let us recall the expression of the proximal operator $\prox_{\lambda f_{\sigma}}$ of PiE \eqref{PiE}, which was systematically explored in \cite{LiuZhouLin} by means of the Lambert W function. The Lambert W function $W(x)$ is a set of solutions of the equation 
$$
x = W (x)e^{W(x)},\mbox{ for any } x\in [-\frac{1}{e},+\infty).
$$
The function $ W(x)$ is  single-valued for  $x\geq 0$ or $x=-\frac{1}{e}$, and is double-valued for
 $-\frac{1}{e}< x< 0$ (see, Fig. \ref{LambertW}). To discriminate between the two branches when $-\frac{1}{e}< x< 0 $, we use the same notation as in  \cite[Section 1.5]{Mezo2022}
 and denote the branch satisfying $W(x)\geq -1$ and $W(x)\leq -1$
 by $W_0(x)$ and $W_{-1}(x)$, respectively. 
 Such a function is a built-in function in Python (\url{https://docs.scipy.org/doc/scipy/reference/generated/scipy.special.lambertw.html}). Lemma \ref{WProperty} later gives their monotonicity. 
  The readers can refer to the recent monograph \cite[Section 1.5]{Mezo2022} on the Lambert W function to learn more details.
 
 \begin{lemma}\cite[Section 1.6]{Mezo2022}\label{WProperty}
  The Lambert W function $W_0(x)$ is strictly increasing on $[-\frac{1}{e},0)$; however, $W_{-1}(x)$ is strictly decreasing on $[-\frac{1}{e},0)$.
 \end{lemma}

The characterizations of the proximal operator of PiE \eqref{PiE} were presented in \cite[Section 2]{LiuZhouLin}, which were split into two cases: $\lambda\le \sigma^2$ and $\lambda>\sigma^2$. For the sake of completeness, we list those characterizations as follows.

\begin{lemma}\label{ProximalTheorem1}
   	Let $\lambda\le \sigma^2$ and $\tau\in \mathbb{R}$. It holds that
   	\begin{align*}
   		\prox_{\lambda f_{\sigma}}(\tau)=\left\{\begin{array}{cl}
   			\{0\}, & \mbox{ if } |\tau|\leq \frac{\lambda}{\sigma},\\
   		 \{\sign(\tau)x_1(\tau)\}, &\mbox{ otherwise},
   		\end{array}\right.
   	\end{align*}
  where
\(
x_1(\tau):=\sigma W_0(-\frac{\lambda}{\sigma^2}e^{-\frac{|\tau|}{\sigma}})+|\tau|.
\)
  \end{lemma}

\begin{lemma}  \label{ProximalTheorem2}
   	  Let $\lambda>\sigma^2$ and $\tau\in \mathbb{R}$. It holds that
   	\begin{align*}
   	\prox_{\lambda f_{\sigma}}(\tau)=\left\{\begin{array}{cl}
   			\{0\}, & \mbox{ if } |\tau|< \sigma(1+\ln\frac{\lambda}{\sigma^2}),\\
   			\sign(\tau)\arg\min\limits_{x=0,x_1(\tau)}\{\widehat{L}(x,\tau)\}, & \mbox{ if } \sigma(1+\ln\frac{\lambda}{\sigma^2})\leq |\tau|\leq \frac{\lambda}{\sigma},\\
   			\{\sign(\tau)x_1(\tau)\}, & \mbox{ otherwise},
   		\end{array}\right.
   	\end{align*}
   where $\hat{L}(x,\tau):=\lambda (1-e^{-\frac{x}{\sigma}})+\frac{1}{2}(x-|\tau|)^2$ and $x_1(\tau)$ is defined as in Lemma  \ref{ProximalTheorem1}.
   \end{lemma}

Lemma \ref{ProximalTheorem2} can be further reduced to the following result, which shows that $\prox_{\lambda f_{\sigma}}(\tau)$ is single-valued except at some point $\bar{\tau}_{\lambda,\sigma}$ depending upon only the $\lambda$ and $\sigma$. This conclusion will be used in the proof of Theorem \ref{MainTh2}.

\begin{lemma}\label{ProximalTheorem3}
Let  $\lambda> \sigma^2$ and $\tau\in \bR$. Then
\begin{equation*}
	\prox_{\lambda f_{\sigma}}(\tau)=\left\{\begin{array}{ll}
\{0\}, & \mbox{ if }|\tau|\leq \bar{\tau}_{\lambda,\sigma},\\
\{0,x_1(\tau)\}, & \mbox{ if }|\tau|= \bar{\tau}_{\lambda,\sigma},\\
\{\sign(\tau)x_1(\tau)\},& \mbox{ otherwise},
\end{array}
\right.
\end{equation*}
where 
$\bar{\tau}_{\lambda,\sigma}=x^*+\frac{\lambda}{\sigma} e^{-\frac{ x^*}{\sigma}}$
with $x^*\in(0,\sqrt{2\lambda})$ being the solution to the equation $\frac{1}{2}+\lambda \frac{(\frac{x}{\sigma}+1)e^{-\frac{x}{\sigma}}-1}{x^2}=0$ on $(0,\infty)$, and $x_1(\tau)$ is defined as in Lemma  \ref{ProximalTheorem1}.
\end{lemma}

Obviously, according to Lemmas \ref{ProximalTheorem2} and \ref{ProximalTheorem3}, the threshold $\bar{\tau}_{\lambda,\sigma}$ satisfies 
$$
\sigma(1+\ln\frac{\lambda}{\sigma^2})\le \bar{\tau}_{\lambda,\sigma}\le \frac{\lambda}{\sigma}.
$$
Those three points will be frequently used when we explore the iteratively reweighted $\ell_1$ algorithm for computing $\prox_{\lambda f_{\sigma}}$ in the next section.

\section{Analysis of IRL1 for computing $\prox_{\lambda f_{\sigma}}$}\label{Section3}

In this section, we will analyze the IRL1 algorithm to compute the following problem:
\begin{align}\label{PIEProx}
\min_{x\in\bR} \Big\{\lambda f_{\sigma}(x)+\frac{1}{2}(x-\tau)^2 \Big\}.
\end{align}
To solve the problem \eqref{PIEProx}, the nonconvex function $f_{\sigma}$ in the IRL1 algorithm is locally approximated by its linear expansion, namely,
 $$
f_{\sigma}(x)\approx f_{\sigma}(x^{(k)})+\frac{1}{\sigma}e^{-\frac{|x^{(k)}|}{\sigma}}(|x|-|x^{(k)}|),
$$
where $x^{(k)}$ denotes the $k$-th iteration. 
With it, the next iteration $x^{(k+1)}$ for a given $\tau$  is computed by
$$
x^{(k+1)}:=\arg\min_{x\in\bR}\Big\{\frac{1}{2}(x-\tau)^2+  \lambda\Big(f_{\sigma}(x^{(k)})+\frac{1}{\sigma}e^{-\frac{|x^{(k)}|}{\sigma}}(|x|-|x^{(k)}|) \Big)\Big\}.
$$
By removing the terms which do not depend on the variable $x$ in the above expression, we obtain
$$
x^{(k+1)}=\arg\min_{x\in\bR}\Big\{\frac{1}{2}(x-\tau)^2+\frac{\lambda}{\sigma}e^{-\frac{|x^{(k)}|}{\sigma}}|x|\Big\}=\prox_{\frac{\lambda}{\sigma}e^{-\frac{|x^{(k)}|}{\sigma}}|\cdot|}(\tau),
$$
that is, 
$$
x^{(k+1)}=\sign(\tau)\Big(|\tau|-\frac{\lambda}{\sigma}e^{-\frac{|x^{(k)}|}{\sigma}}\Big)_+,
$$
where $(t)_+:=\max\{0,t\}$.

It is sufficient to restrict our discussion on $\tau>0$ as $\prox_{\lambda f_{\sigma}}(\tau)$ is symmetric about the origin  \cite[Lemma 2.1]{LiuZhouLin} and $\prox_{\lambda f_{\sigma}}(0)=\{0\}$. 
To be more precise,  the IRL1 algorithm for PiE with $\tau>0$ is described in Algorithm 1.

\begin{algorithm}
\caption{Iteratively Reweighted $\ell_1$ Algorithm (IRL1)}\label{Algorithm}
{\bf Input} Fix $\lambda>0$ and $\sigma>0$. Given $x^{(0)}\ge0$ and $\tau>0$.
\begin{itemize}
\item[] {\bf for} $k=0,1,\dots$ {\bf do}
\begin{equation}\label{reweightedl1}
x^{(k+1)}=\Big(\tau-\frac{\lambda}{\sigma}e^{-\frac{x^{(k)}}{\sigma}}\Big)_+
\end{equation}
\item[] {\bf end for}
\end{itemize}
{\bf Output} $x^{(\infty)}$
\end{algorithm}

Denote $F(x)\!:=\!\lambda f_{\sigma}(x)+\frac{1}{2}(x-\tau)^2$. We call $x$ a critical point of the function $F$, if  $0\!\in\! \partial F(x)$ is satisfied, where $\partial F(x)$ denotes the subdifferential of $F$ at $x$ \cite[Definition 8.3]{Rockafellar2009}.
Ochs et al \cite{Ochs2015} pointed out that the sequence
$\{x^{(k)}\}$ generated by Algorithm \ref{Algorithm}
converges to a critical point of the function $F$.
We go
one step further than the previous result and show that not only the sequence $\{x^{(k)}
\}$ is convergent, but also its limit $x^{(\infty)}$ depends on the initialization $x^{(0)}$ and the relationship of $\tau$ with
the parameters $\lambda$ and $\sigma$.
The convergence behavior of \eqref{reweightedl1} is described by Lemmas \ref{Exik0Lem}--\ref{ComplexLemCase2} in Section \ref{Section4}.
 This is then compared to the true solution set $\prox_{\lambda f_{\sigma}}(\tau)$ in Theorems \ref{MainTh1} and \ref{MainTh2}. In particular,
we identify the intervals where \eqref{reweightedl1} will not achieve the true solution. These intervals are
explicitly determined in terms of the initial  $x^{(0)}$ and parameters $\lambda$ and $\sigma$.

Notice that $x^{(\infty)}$ satisfying the equation $x=(\tau-\frac{\lambda}{\sigma}e^{-\frac{x}{\sigma}})_+$ by \eqref{reweightedl1}. To further investigate properties of $x^{(\infty)}$,  for given $\tau\!\in\! \mathbb{R}$ we define a function $\phi:\mathbb{R}\to \mathbb{R}$  with
\begin{equation}\label{phiDef}
\phi(x):=\tau-x-\frac{\lambda}{\sigma}e^{-\frac{x}{\sigma}}, {\text{ for any }} x\in \mathbb{R},
\end{equation}
 and its main properties used later are listed in the following Lemma.
\begin{lemma}\label{phiLem} Let $\phi$ be defined by \eqref{phiDef}. Write
  \(
  x_2(\tau):=\sigma W_{-1}(-\frac{\lambda}{\sigma^2}e^{-\frac{\tau}{\sigma}})+\tau
  \), and $x_1(\tau)$ is defined as in Lemma  \ref{ProximalTheorem1}.
Then, the following statements hold.
\begin{itemize}
    \item [{\bf (i)}] The function  $\phi$
is strictly increasing on $(-\infty,\sigma\ln\frac{\lambda}{\sigma^2}]$ and strictly decreasing on $(\sigma\ln\frac{\lambda}{\sigma^2},+\infty)$.
Moreover, $\phi(x)\!\le\! \phi(\sigma\ln\frac{\lambda}{\sigma^2})\!=\!\tau\!-\!\sigma (1\!+\!\ln \frac{\lambda}{\sigma^2})$ for any $x\!\in\! \mathbb{R}$.
 
    \item [{\bf (ii)}] If $\tau\!\in\! (\sigma(1\!+\!\ln\frac{\lambda}{\sigma^2}),\frac{\lambda}{\sigma})$, the equation $\phi(x)\!=\!0$  has two solutions $x_1(\tau)$ and $x_2(\tau)$ with
  \begin{align}
      \left\{\begin{array}{cl}
      0<x_2(\tau)< \sigma\ln{\frac{\lambda}{\sigma^2}}<x_1(\tau), & {\rm if\;} \lambda> \sigma^2, \\
      x_2(\tau)<\sigma\ln{\frac{\lambda}{\sigma^2}}<x_1(\tau)<0, & {\rm if\;} \lambda\leq \sigma^2. 
      \end{array}\right.
  \end{align}
  
   \item [{\bf (iii)}]  If  $\tau= \sigma(1\!+\!\ln\frac{\lambda}{\sigma^2})$, the equation $\phi(x)=0$  has a unique solution, that is,  $x_1(\tau)=x_2(\tau)=\sigma\ln{\frac{\lambda}{\sigma^2}}$.
  
  \item [{\bf (iv)}] If $\tau>\frac{\lambda}{\sigma}$, the equation $\phi(x)=0$  has two solutions $x_1(\tau)$ and $x_2(\tau)$ satisfying $x_2(\tau)<0<x_1(\tau)$.

  \item [{\bf (v)}] If $\tau=\frac{\lambda}{\sigma}$,
  the equation $\phi(x)=0$  has two solutions $x_1(\tau)$ and $x_2(\tau)$ with
  \begin{align}
      \left\{\begin{array}{cl}
      0=x_2(\tau)<\sigma\ln{\frac{\lambda}{\sigma^2}}<x_1(\tau), & {\rm if\;} \lambda>\sigma^2, \\
      x_2(\tau)<\sigma\ln{\frac{\lambda}{\sigma^2}}<x_1(\tau) =0,& {\rm if\;} \lambda<\sigma^2, \\
      x_1(\tau)=x_2(\tau)=0, & {\rm if\;} \lambda=\sigma^2. 
      \end{array}\right.
  \end{align}
  
\end{itemize}
\end{lemma}
\begin{proof}
  After simple calculation, 
 $\phi'(x)=\frac{\lambda}{\sigma^2}e^{-\frac{x}{\sigma}}-1$, $\phi''(x)=-\frac{\lambda}{\sigma^3}e^{-\frac{x}{\sigma}}$. Clearly, 
 the statement (i) holds.
 The equation $\phi(x)=0$ is equivalent to   $x\!=\!\tau-\frac{\lambda}{\sigma}e^{-\frac{x}{\sigma}}$, namely, 
 \begin{equation}\label{phi(x)=0}
\frac{x-\tau}{\sigma}e^{\frac{x-\tau}{\sigma}}=-\frac{\lambda}{\sigma^2}e^{-\frac{\tau}{\sigma}}.
\end{equation}
If $\tau\!>\!\sigma(1\!+\!\ln\frac{\lambda}{\sigma^2})$, $-\frac{\lambda}{\sigma^2}e^{-\frac{\tau}{\sigma}}\!\in\! (-\frac{1}{e},0)$. By definition of Lambert W function and the equation \eqref{phi(x)=0}, the equation $\phi(x)=0$  has two solutions $x_1(\tau)$ and $x_2(\tau)$.
Together with (i) 
and the fact $\phi(0)\!=\!\tau-\frac{\lambda}{\sigma}$,  we know that the statements (ii) and (iv) hold.
 When $\tau\!=\!\sigma(1\!+\!\ln\frac{\lambda}{\sigma^2})$,  $-\frac{\lambda}{\sigma^2}e^{-\frac{\tau}{\sigma}}\!=\!-\frac{1}{e}$.
Hence, we obtain
\[
W_{-1}(-\frac{\lambda}{\sigma^2}e^{-\frac{\tau}{\sigma}})\!=\!W_{0}(-\frac{\lambda}{\sigma^2}e^{-\frac{\tau}{\sigma}})\!=\!W(-\frac{1}{e})=-1,
\]
 which implies $x_1(\tau)\!=\!x_2(\tau)\!
 =\!\tau-\sigma\!=\!\sigma\ln{\frac{\lambda}{\sigma^2}}$. The statement (iii) holds. In the following, we will argue the  statement (v).
 Notice
 $\tau\!=\!\frac{\lambda}{\sigma}\!>\! \sigma(1\!+\!\ln\frac{\lambda}{\sigma^2})$, $-\frac{\lambda}{\sigma^2}e^{-\frac{\tau}{\sigma}}\!\in\! (-\frac{1}{e},0)$.
 So, the equation $\phi(x)\!=\!0$  has  solutions $x_1(\tau)$ and $x_2(\tau)$.  From (i), it follows
 \begin{align}\label{phiTemp1}
 x_2(\tau)<\sigma \ln{\frac{\lambda}{\sigma^2}}<x_1(\tau).
 \end{align}
 We will proceed in two cases.
 
 {\bf Case 1: $\lambda\neq \sigma^2$.}
 If $\lambda>\sigma^2$, then 
 $-\frac{\lambda}{\sigma^2}<-1$.
 With
 $-\frac{\lambda}{\sigma^2}e^{-\frac{\tau}{\sigma}}\!\in\! (-\frac{1}{e},0)$, we know that
 \(
 W_{-1}(-\frac{\lambda}{\sigma^2}e^{-\frac{\lambda}{\sigma^2}})=-\frac{\lambda}{\sigma^2}.
 \)
 Together with $\tau=\frac{\lambda}{\sigma}$,  yielding
 \begin{align*}
x_2(\tau)=\sigma W_{-1}(-\frac{\lambda}{\sigma^2}e^{-\frac{\tau}{\sigma}})+\tau
=\tau+\sigma W_{-1}(-\frac{\lambda}{\sigma^2}e^{-\frac{\lambda}{\sigma^2}})
=\tau-\frac{\lambda}{\sigma}=0.
\end{align*}
Again from \eqref{phiTemp1}, it follows that
$0=x_2(\tau)<\sigma \ln{\frac{\lambda}{\sigma^2}}<x_1(\tau)$.

If $\lambda<\sigma^2$, then 
 $-\frac{\lambda}{\sigma^2}>-1$.
 With
 $-\frac{\lambda}{\sigma^2}e^{-\frac{\tau}{\sigma}}\!\in\! (-\frac{1}{e},0)$, we know that
 \(
 W_{0}(-\frac{\lambda}{\sigma^2}e^{-\frac{\lambda}{\sigma^2}})=-\frac{\lambda}{\sigma^2}.
 \)
 Together with $\tau=\frac{\lambda}{\sigma}$, yielding
 \begin{align*}
x_1(\tau)=\sigma W_{0}(-\frac{\lambda}{\sigma^2}e^{-\frac{\tau}{\sigma}})+\tau
=\tau+\sigma W_0(-\frac{\lambda}{\sigma^2}e^{-\frac{\lambda}{\sigma^2}})
=\tau-\frac{\lambda}{\sigma}=0.
\end{align*}
Again from \eqref{phiTemp1}, it follows  that
$x_2(\tau)<\sigma \ln{\frac{\lambda}{\sigma^2}}<x_1(\tau)=0$.

{\bf Case 2:} $\lambda\!=\!\sigma^2$. Now $-\frac{\lambda}{\sigma^2}e^{-\frac{\tau}{\sigma}}\!=\!-\frac{1}{e}$ by $\tau=\frac{\lambda}{\sigma}$.
Hence,
\[
W_{-1}(-\frac{\lambda}{\sigma^2}e^{-\frac{\tau}{\sigma}})\!=\!W_{0}(-\frac{\lambda}{\sigma^2}e^{-\frac{\tau}{\sigma}})\!=\!W(-\frac{1}{e})=-1,
\]
 which implies   $\phi(x)\!=\!0$  has a unique solution $x_1(\tau)\!=\!x_2(\tau)\!=\!0$
 by \eqref{phi(x)=0}.
\end{proof}

\begin{proposition}\label{xinftyProp}
  Given $\tau>0$ and an initial value $x^{(0)}\ge0$.
  Suppose that the sequence $\{x^{(k)}\}$ generated by Algorithm \ref{Algorithm} converges to  $x^{(\infty)}$. Then, the following statements hold.
 \begin{itemize}
     \item [{\bf(i)}]  $x^{(\infty)}=0$ implies that $\tau\leq \frac{\lambda}{\sigma}$.
     
     \item [{\bf(ii)}] If  $\tau> \frac{\lambda}{\sigma}$, $x^{(\infty)}=\sigma W_0(-\frac{\lambda}{\sigma^2}e^{-\frac{\tau}{\sigma}})+\tau$.
 \end{itemize}
\end{proposition}
 \begin{proof}
 By the continuity of the function $(\cdot)_+$, $x^{(k)}\to x^\infty$ and the equation \eqref{reweightedl1} for each $k$, we know that
 \begin{align}\label{x_inftyEQ}
x^{(\infty)}=(\tau-\frac{\lambda}{\sigma}e^{-\frac{x^{(\infty)}}{\sigma}})_+.
\end{align}
If $x^{(\infty)}\!=\!0$, then
$(\tau-\frac{\lambda}{\sigma})_+\!=\!0$ from \eqref{x_inftyEQ}, which implies that
 $\tau\!\leq \!\frac{\lambda}{\sigma}$.
 Hence, the statement (i) holds. 
 If $\tau\!>\! \frac{\lambda}{\sigma}$, then $x^{(\infty)}\!>\!0$ from (i), and 
 $x^{(\infty)}\!=\!\tau-\frac{\lambda}{\sigma}e^{-\frac{x^{(\infty)}}{\sigma}}$ from \eqref{x_inftyEQ}, namely, $\phi(x^{(\infty)})\!=\!0$, where $\phi$ defined by \eqref{phiDef}. So, 
 $x^{(\infty)}\!=\!\sigma W_0(-\frac{\lambda}{\sigma^2}e^{-\frac{\tau}{\sigma}})+\tau$ from Lemma \ref{phiLem} (iv).
 \end{proof}

\subsection{Comparing IRL1 solution with $\prox_{\lambda f_{\sigma}}$}

In this subsection, 
 we will identify when the limit $x^{(\infty)}$ of the sequence $\{x^{(k)}\}$ belongs  or not belongs to the set $\prox_{\lambda f_{\sigma}}(\tau)$.  We recall in Lemmas \ref{ProximalTheorem1} and \ref{ProximalTheorem3} that the set $\prox_{\lambda f_{\sigma}}(\tau)$ has a unique element except for $|\tau|=\bar{\tau}_{\lambda,\sigma}$ with $\lambda>\sigma^2$.

The following two theorems summarize our main results. Our results for PiE are mainly inspired by the ideas presented in \cite[section 4]{Prater2022} for the iteratively reweighted algorithm for computing the proximal operator of the log-sum penalty.  The proofs as well as relevant technical lemmas are given in Section \ref{Section4}. 
From now on, we say that a sequence $\{x^{(k)}\}$ is converging to $\prox_{\lambda f_{\sigma}}(\tau)$ provided that the limit of $\{x^{(k)}\}$ belongs to the set  $\prox_{\lambda f_{\sigma}}(\tau)$.

\begin{theorem}\label{MainTh1}
Given $\tau\!>\! 0$ and an initial value $x^{(0)}\!\geq \!0$.
Let $\lambda\le \sigma^2$. Then the sequence $\{x^{(k)}\}$ generated by Algorithm \ref{Algorithm}  converges to $\prox_{\lambda f_{\sigma}}(\tau)$.
\end{theorem}

\begin{figure}[!tbp]
\begin{center}
\includegraphics[scale=0.025]{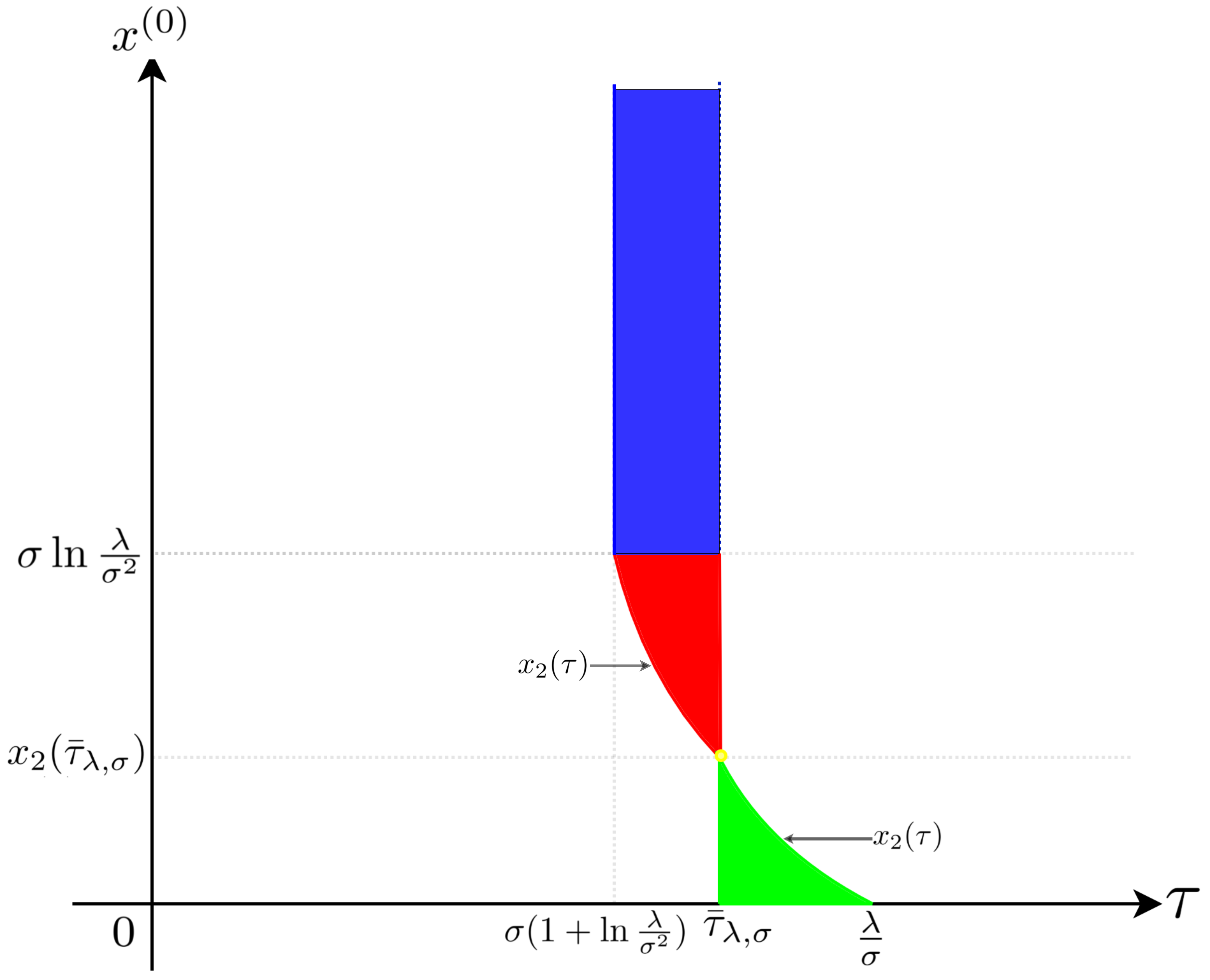}
\caption{Illustration of Theorem \ref{MainTh2} with $\tau>0$. The specific regions where the IRL1 solution differs from the proximal operator of PiE in (ii)-(v) are marked in blue, red, yellow, and green, respectively.}
\label{IRL1}
\end{center}
\end{figure}

If $\lambda\!>\!\sigma^2$, we see that $\{x^{(k)}\}$ generated by Algorithm \ref{Algorithm} may not always converge to  $\prox_{\lambda f_{\sigma}}$ for some given $x^{(0)}\!\ge\!0$. The regions where the algorithm fails depend on the threshold $\bar{\tau}_{\lambda,\sigma}$ given as in Lemma \ref{ProximalTheorem3} and $x_2(\tau)$ defined in Lemma \ref{phiLem}, as shown in Fig. \ref{IRL1}. The value $\bar{\tau}_{\lambda,\sigma}$ can be computed by the bisection method. 
Notice that $x_2(\tau)$ is  strictly decreasing on $[\sigma(1+\ln\frac{\lambda}{\sigma^2}),\frac{\lambda}{\sigma}]$, by Lemma \ref{WProperty},   we denote the inverse function of $x_2(\tau)$  by $x_2^{-1}(\tau)$ for each $\tau\in [\sigma(1+\ln\frac{\lambda}{\sigma^2}),\frac{\lambda}{\sigma}]$.

\begin{theorem}\label{MainTh2}  Given $\tau\!>\! 0$ and an initial value $x^{(0)}\!\geq \!0$.
Let $\lambda> \sigma^2$, $\bar{\tau}_{\lambda,\sigma}$  be defined as in Lemma \ref{ProximalTheorem3},
$x_i(\tau)(i=1,2)$ be defined as in Lemma \ref{phiLem} and the sequence $\{x^{(k)}\}$ be generated by Algorithm \ref{Algorithm}. Then the following statements hold.
\begin{itemize}

\item[{\bf(i)}] The sequence $\{x^{(k)}\}$ converges
to $\prox_{\lambda f_{\sigma}}(\tau)$ for any $\tau\!\in\! (0,\sigma(1+\ln \frac{\lambda}{\sigma^2}))\!\cup\! (\frac{\lambda}{\sigma},+\infty)$.

\item[{\bf(ii)}] If $x^{(0)}\ge \sigma\ln\frac{\lambda}{\sigma^2}$, 
$\{x^{(k)}\}$ converges to $x_1(\tau)$ for any $\tau\in [\sigma(1+\ln \frac{\lambda}{\sigma^2}),\frac{\lambda}{\sigma}]$. Consequently, $\{x^{(k)}\}$ converges
to $\prox_{\lambda f_{\sigma}}(\tau)$ for any $\tau\in [\bar{\tau}_{\lambda,\sigma},\frac{\lambda}{\sigma}]$, however 
 $\{x^{(k)}\}$ does not converge to  $\prox_{\lambda f_{\sigma}}(\tau)$ for any $\tau\in [\sigma(1+\ln\frac{\lambda}{\sigma^2}),\bar{\tau}_{\lambda,\sigma})$.
 
\item[{\bf(iii)}] If $x_2(\bar{\tau}_{\lambda,\sigma})\!<\!x^{(0)}\!<\! \sigma\ln\frac{\lambda}{\sigma^2}$, the sequence $\{x^{(k)}\}$ converges
to $\prox_{\lambda f_{\sigma}}(\tau)$ for any $\tau\in [\sigma(1+\ln \frac{\lambda}{\sigma^2}),x_2^{-1}(x^{(0)}))\cup [\bar{\tau}_{\lambda,\sigma},\frac{\lambda}{\sigma}]$, but the sequence
$\{x^{(k)}\}$ does not converge to  $\prox_{\lambda f_{\sigma}}(\tau)$ for any $\tau\in [x_2^{-1}(x^{(0)}), \bar{\tau}_{\lambda,\sigma})$.

\item[{\bf(iv)}] If $x^{(0)}=x_2(\bar{\tau}_{\lambda,\sigma})$, the sequence $\{x^{(k)}\}$ converges to $\prox_{\lambda f_{\sigma}}(\tau)$ for any $\tau\in [\sigma(1+\ln \frac{\lambda}{\sigma^2}), \bar{\tau}_{\lambda,\sigma})
\cup(  \bar{\tau}_{\lambda,\sigma},\frac{\lambda}{\sigma}]$, however the sequence $\{x^{(k)}\}$ does not converges to $\prox_{\lambda f_{\sigma}}(\tau)$ when 
$\tau= \bar{\tau}_{\lambda,\sigma}$.

\item[{\bf(v)}] If $0\!\le\! x^{(0)}\!<\!x_2(\bar{\tau}_{\lambda,\sigma})$, the sequence $\{x^{(k)}\}$ converges
to $\prox_{\lambda f_{\sigma}}(\tau)$ for any $\tau\in [\sigma(1+\ln \frac{\lambda}{\sigma^2}),\bar{\tau}_{\lambda,\sigma}]\cup (x_2^{-1}(x^{(0)}),\frac{\lambda}{\sigma}]$, but the sequence
$\{x^{(k)}\}$ does not  converge to  $\prox_{\lambda f_{\sigma}}(\tau)$ for any $\tau\in ( \bar{\tau}_{\lambda,\sigma},x_2^{-1}(x^{(0)})]$.
\end{itemize}
\end{theorem}

\begin{remark} The initial value for ILR1 is usually and simply set to be $1$ \cite[Subsection 2.2]{Candes2008} for compressed sensing, to be a random feasible value for support vector machines \cite[Subsection 2.1]{Bradley1998b}, and the identity matrix for a low-rank matrix completion problem \cite[Algorithm1]{Zhou2022}. By Theorem \ref{MainTh2}, the above choice may result in the deviation between the IRL1 solution and the proximal operator of PiE.
\end{remark}

\begin{figure}[!tbp]
\begin{center}
\includegraphics[scale=0.5]{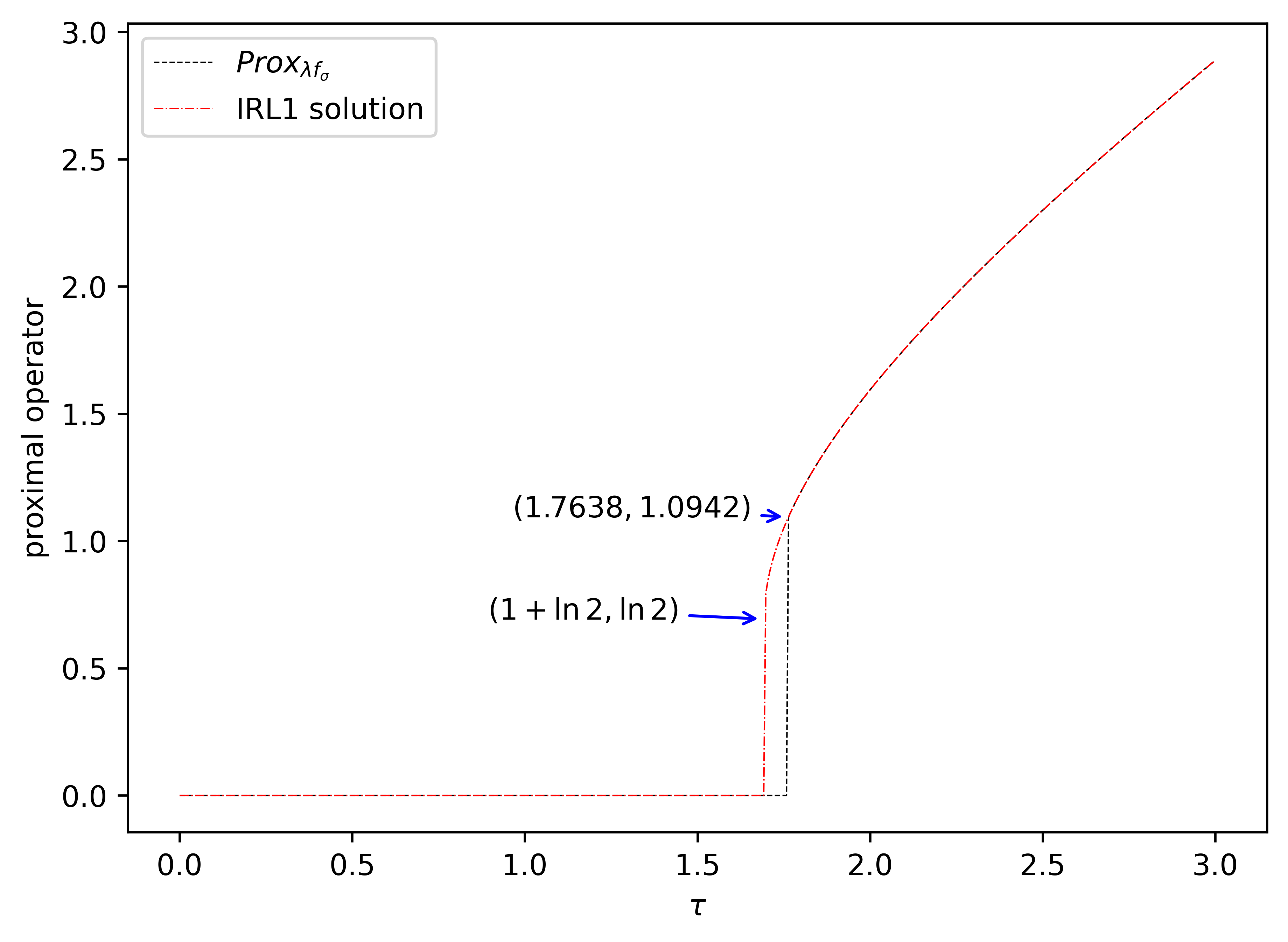}
\caption{The proximal operator $\prox_{\lambda f_{\sigma}}$ and the IRL1 solution in Theorem \ref{MainTh2} with $\lambda=2$, $\sigma=1$, and the initial value $x^{(0)}=1$.}
\label{Trueproximal}
\end{center}
\end{figure}

Fig. \ref{IRL1} illustrates the results (i)-(v) in Theorem \ref{MainTh2} with $\tau\!>\!0$.  Only when $\tau$ lies in a subset of the interval $[\sigma(1+\ln \frac{\lambda}{\sigma^2}),\frac{\lambda}{\sigma}]$ does the deviation occur. The colored regions indicate where the IRL1 solution differs from the proximal operator of PiE.
For example, let $\lambda=2$ and $\sigma=1$. Then $\sigma(1+\ln\frac{\lambda}{\sigma^2})=1+\ln 2$, $\bar{\tau}_{\lambda,\sigma}=1.7638$, $\frac{\lambda}{\sigma}=2$, $\sigma\ln\frac{\lambda}{\sigma^2}=\ln 2$,  $x_2(\bar{\tau}_{\lambda,\sigma})=0.3393$, and $x_1(\bar{\tau}_{\lambda,\sigma})=1.094$.  
In this case, given an initial value $x^{(0)}=1>\sigma\ln\frac{\lambda}{\sigma^2}$, the IRL1 solution (red dashdot) and the true proximal operator (black dashed) are illustrated in Fig. \ref{Trueproximal}, which corresponds to the case of Theorem \ref{MainTh2} (ii). Clearly, the IRL1 solution disagrees with the true proximal operator for any given $\tau\in[1+\ln 2,1.7638)$. 

\subsection{Choices of initial values}

We are devoted to adaptively selecting an initial value in a simple way to guarantee the fast convergence of the IRL1 solution to $\prox_{\lambda f_{\sigma}}(\tau)$ for all $\tau>0$.  The discussion will be divided into two cases: $\lambda\le\sigma^2$ and $\lambda>\sigma^2$.

\begin{theorem}\label{Theorem3} 
Given $\tau>0$.
Suppose that $\lambda\le\sigma^2$ and the sequence $\{x^{(k)}\}$ is generated by Algorithm \ref{Algorithm} with the initial value $x^{(0)}$ in  Algorithm \ref{Algorithm} given as 
\begin{equation}\label{initialsoft}
x^{(0)}:=\left\{\begin{array}{ll}
0,& \mbox{ if } \tau\le \frac{\lambda}{\sigma},\\
\tau, &\mbox{ otherwise}.
\end{array}
\right.
\end{equation}
Then, the following statements hold.
\begin{itemize}
    \item [{\bf(i)} ] If $0<\tau\le \frac{\lambda}{\sigma}$, then $x^{(k)}=0$  for each $k\in\bN$.
    
    \item [{\bf(ii)} ] If $\tau>\frac{\lambda}{\sigma}$, it holds that
\begin{equation}\label{exponential}
\Big(\frac{\lambda}{\sigma^2}e^{-\frac{\tau}{\sigma}}\Big)^k (\tau-x_1(\tau))\!<\! x^{(k)}-x_1(\tau)\!<\!\Big(\frac{\lambda}{\sigma^2}e^{-\frac{x_1(\tau)}{\sigma}}\Big)^k (\tau-x_1(\tau)),
\end{equation}
where $x_1(\tau)$ is defined as in Lemma  \ref{ProximalTheorem1}.
\end{itemize}
\end{theorem}
\begin{proof}
   If $0\!<\!\tau\le \frac{\lambda}{\sigma}$, $x^{(0)}\!=\!0$, the statement (i) is trivial by Lemma \ref{Exik0Lem} (i). 
   Now suppose $\tau\!>\! \frac{\lambda}{\sigma}$. Then $x^{(0)}\!=\!\tau>0$.  Therefore,
    by Lemma \ref{x0geq0Lem},
  $x^{(k)}>0$ for each 
$k\in \mathbb{N}$ and
 $\{x^{(k)}\}$  converges to $x_1(\tau)$.
 Notice that
 $\tau=x_1(\tau)+\frac{\lambda}{\sigma}e^{-\frac{x_1(\tau)}{\sigma}}$ 
 by Lemma \ref{phiLem}(iv).
 With \eqref{reweightedl1}
 and the Lagrange mean value theorem, we arrive at
\begin{equation}\label{xtau1}
x^{(k+1)}-x_1(\tau)=\frac{\lambda}{\sigma}\Big(e^{-\frac{x_1(\tau)}{\sigma}}-e^{-\frac{x^{(k)}}{\sigma}}\Big)=\frac{\lambda}{\sigma^2}e^{-\frac{\xi}{\sigma}}(x^{(k)}-x_1(\tau)),
\end{equation}
for some $\xi\in (x_1(\tau),x^{(k)})\subseteq (x_1(\tau),\tau)$. Note that $e^{-\frac{t}{\sigma}}$ is strictly decreasing for any $t>0$.  By \eqref{xtau1}, it follows that
$$
\frac{\lambda}{\sigma^2}e^{-\frac{\tau}{\sigma}}(x^{(k)}-x_1(\tau))\!<\!x^{(k+1)}-x_1(\tau)\!<\!\frac{\lambda}{\sigma^2}e^{-\frac{x_1(\tau)}{\sigma}}(x^{(k)}-x_1(\tau)).
$$
Repeating the process, which yields \eqref{exponential}. 
The proof is hence complete.
\end{proof}

  \begin{remark}
  A sequence $\{y^{(k)}\}\subseteq \bR$ is said to converge Q-linearly to a point $\bar{y}$ if there exists $c>0$ such that $\lim_{k\to+\infty}|y^{(k+1)}-\bar{y}|/|y^{(k)}-\bar{y}|=c$. 
    The equation \eqref{xtau1} in Theorem \ref{Theorem3}  shows that  the approximate error $x^{(k)}-\prox_{\lambda f_{\sigma}}$ converges Q-linearly to $x_1(\tau)$ for each $\tau\!>\!\frac{\lambda}{\sigma}$ when $\lambda\leq \sigma^2$.
  \end{remark}

\begin{figure}
\centering
\subfloat[]{\includegraphics[width=0.45\textwidth]{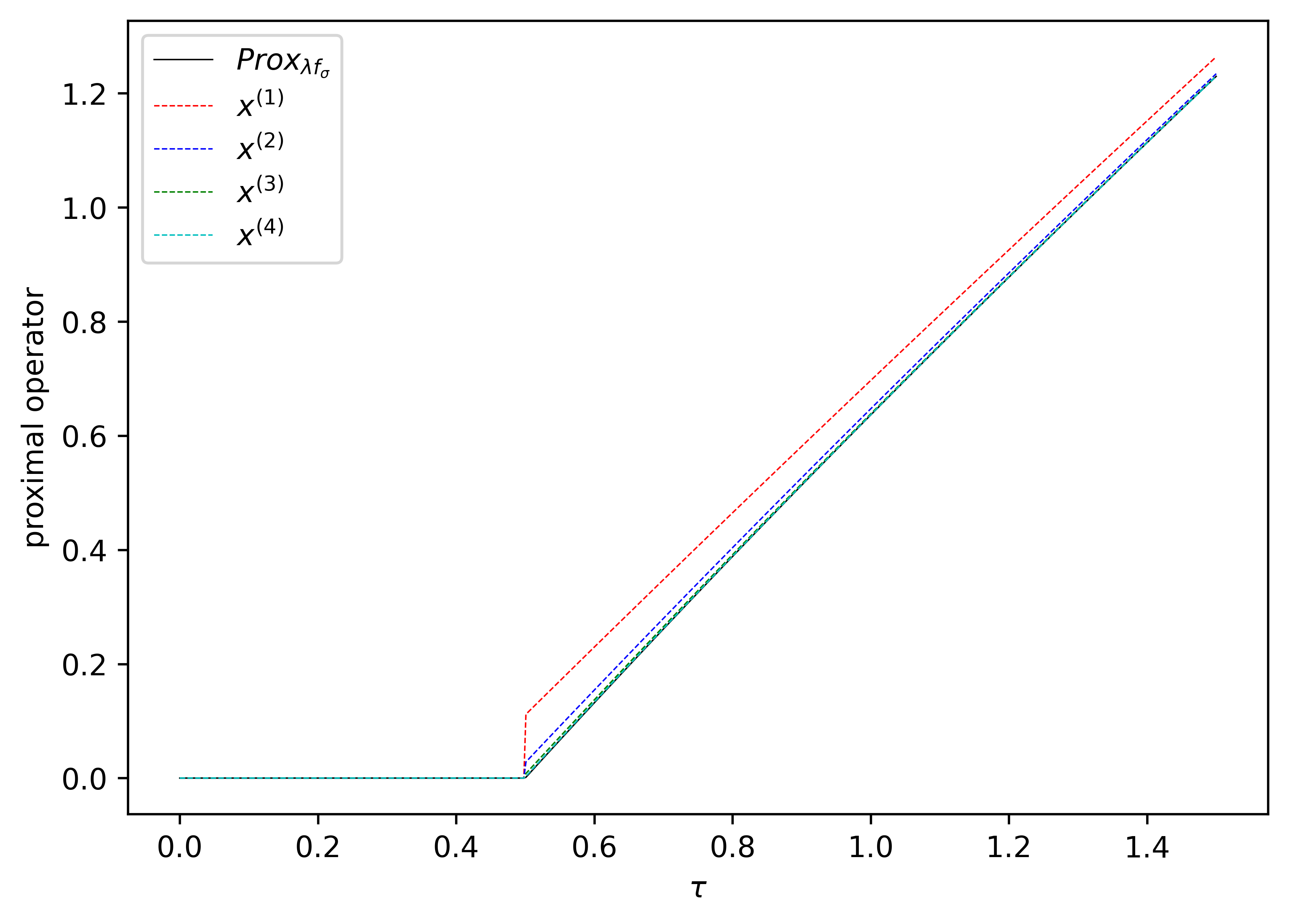}}
\subfloat[]{\includegraphics[width=0.45\textwidth]{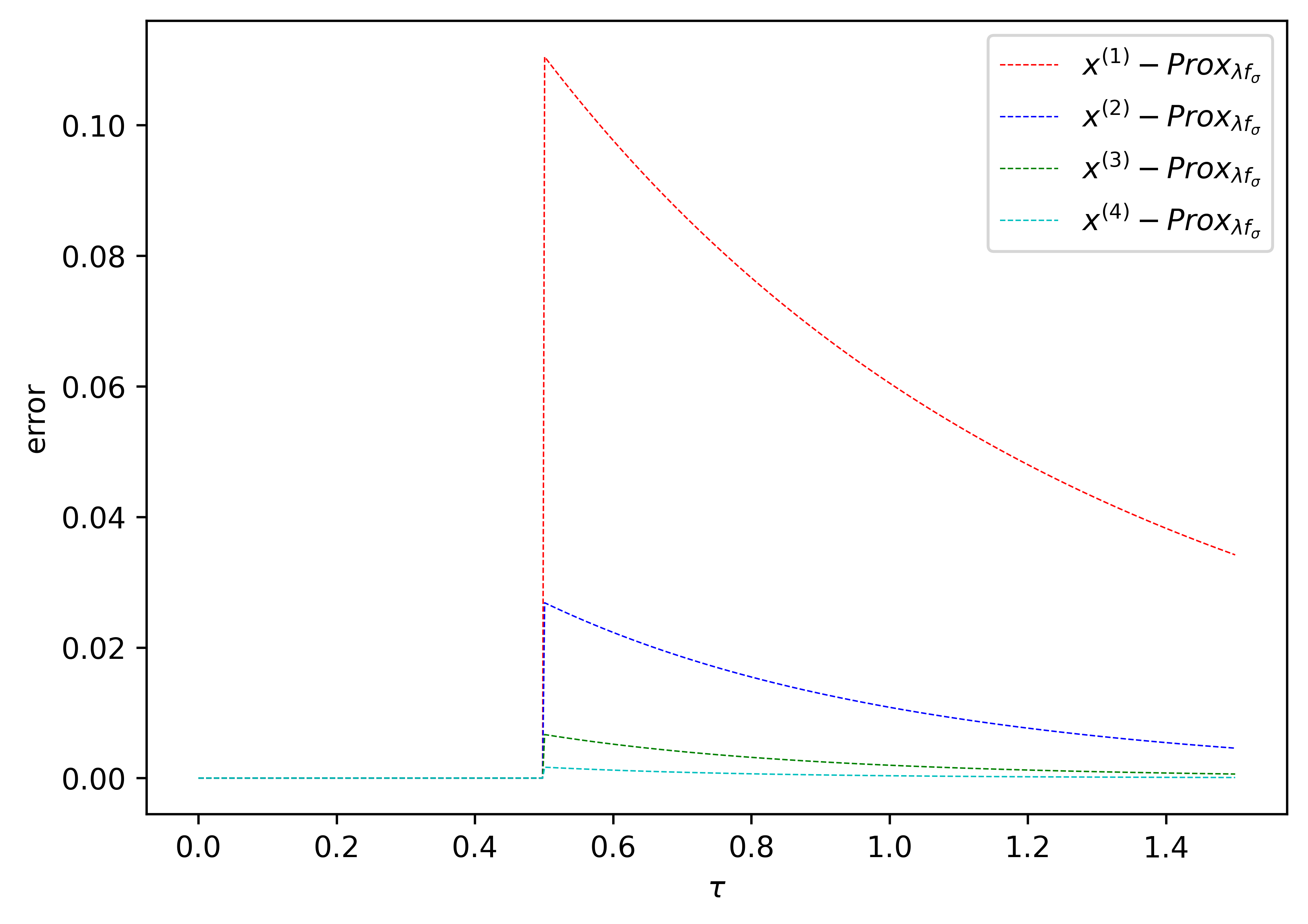}}
\caption{Illustration of Theorem \ref{Theorem3} with $\lambda=1$ and $\sigma=2$ for any $\tau\in(0,1.5]$. (a) The true proximal operator and IRL1 solution $x^{(k)}$ with $k=1,2,3,4$; (b) the error function.}
\label{PiE_prox_soft}
\end{figure}

For \eqref{exponential}, 
$\frac{\lambda}{\sigma^2}e^{-\frac{x_1(\tau)}{\sigma}}< \frac{\lambda}{\sigma^2}e^{0}\le 1$ by Lemma \ref{phiLem}  (iv). Moreover,
$x_1(\tau)$ is increasing on $(\frac{\lambda}{\sigma},+\infty)$ by Lemma \ref{WProperty} and $x_1(\tau)\to 0$ as $\tau\to \frac{\lambda}{\sigma}$.
Let $\lambda=1$ and $\sigma=2$ in Theorem \ref{Theorem3}.
Fix $k\in\{1,2,3,4\}$, $x^{(k)}$, $\prox_{\lambda f_{\sigma}}$, and the corresponding error function $x^{(k)}-\prox_{\lambda f_{\sigma}}$ for any $\tau>0$ are illustrated in Fig. \ref{PiE_prox_soft}. 

Next, we will consider the case that $\lambda>\sigma^2$.
As a direct sequence of Theorem \ref{MainTh2}, if we fix the initial value $x^{(0)}\ge0$ (for example, $x^{(0)}=\tau$ or $x^{(0)}=0$) for all $\tau>0$, then
$\{x^{(k)}\}$ generated by Algorithm \ref{Algorithm} fails to converge to the solution of $\prox_{\lambda f_{\sigma}}(\tau)$ for at least one $\tau>0$. 
To solve this, 
 the initial value $x^{(0)}\ge0$ will be chose  depending on $\bar{\tau}_{\lambda,\sigma}$. 
A simple choice for $x^{(0)}$ is suggested below.

\begin{theorem}\label{Theorem4}
Given $\tau\!>\!0$. Let  $\lambda\!>\!\sigma^2$ and the initial value $x^{(0)}$ is given by
\begin{equation}\label{initialhard}
x^{(0)}:=\left\{\begin{array}{ll}
0,& \mbox{ if } \tau\le \bar{\tau}_{\lambda,\sigma},\\
\tau, &\mbox{ otherwise}.
\end{array}
\right.
\end{equation}
Then, the following statements hold.
\begin{itemize}
    \item [{\bf (i) }] The sequence $\{x^{(k)}\}$ generated by Algorithm \ref{Algorithm} converges to the solution of $\prox_{\lambda f_{\sigma}}(\tau)$ for any $\tau>0$.
    
    \item [{\bf(ii)} ] If $0<\tau\le \bar{\tau}_{\lambda,\sigma}$, then $x^{(k)}=0$  for each $k$.
    
    \item [{\bf(iii)} ] If $\tau>\bar{\tau}_{\lambda,\sigma}$, it holds that
\begin{equation*}\label{exponential2}
\Big(\frac{\lambda}{\sigma^2}e^{-\frac{\tau}{\sigma}}\Big)^k (\tau-x_1(\tau))\!<\! x^{(k)}-x_1(\tau)\!<\!\Big(\frac{\lambda}{\sigma^2}e^{-\frac{x_1(\tau)}{\sigma}}\Big)^k (\tau-x_1(\tau)),
\end{equation*}
where $x_1(\tau)$ is defined as in Lemma  \ref{ProximalTheorem1}.
\end{itemize}

\end{theorem}
\begin{proof} 
Suppose $\tau\le \bar{\tau}_{\lambda,\sigma}$. Then
 $x^{(0)}=0$. By Theorem \ref{MainTh2} (i) and (v), then $\{x^{(k)}\}$ converges to  $\prox_{\lambda f_{\sigma}}(\tau)$ for any $\tau\le \bar{\tau}_{\lambda,\sigma}$. If $\tau>\bar{\tau}_{\lambda,\sigma}$, then $\tau>\sigma(1+\ln\frac{\lambda}{\sigma^2})$ and further $x^{(0)}=\tau>\sigma\ln\frac{\lambda}{\sigma^2}$. By Theorem \ref{MainTh2} (i) and (ii), $\{x^{(k)}\}$ converges to  $\prox_{\lambda f_{\sigma}}(\tau)$ for any $\tau\!>\!\bar{\tau}_{\lambda,\sigma}$.
 So,  the statement (i) holds. 
 The statement (ii) is from Lemma \ref{Exik0Lem} (i)
 and the fact 
 $\bar{\tau}_{\lambda,\sigma}\!\leq \!\frac{\lambda}{\sigma}$. Now suppose that $\tau>\bar{\tau}_{\lambda,\sigma}$.  Then $x^{(0)}=\tau$. Notice $\phi(\tau)<0$ where $\phi$ is defined in \eqref{phiDef}. Then $x_1(\tau)<\tau$ by Lemma \ref{phiLem} and $\lambda>\sigma^2$. Associating  the proof of Lemma \ref{ComplexLemCase2} (iii) when $x^{(0)}>x_1(\tau)$ with Lemma \ref{x0geq0Lem}, 
   we know that
  \(
 x^{(k+1)}\!=\!\tau-\frac{\lambda}{\sigma}e^{-\frac{x^{(k)}}{\sigma}}
 \)
 and $x^{(k)}\!>\! x^{(k+1)}\!>\!x_1(\tau)$ for each $k$.
 The rest proof is similar to the last part of Theorem \ref{Theorem3}. We omit it.
\end{proof}

\begin{figure}
\centering
\subfloat[]{\includegraphics[width=0.45\textwidth]{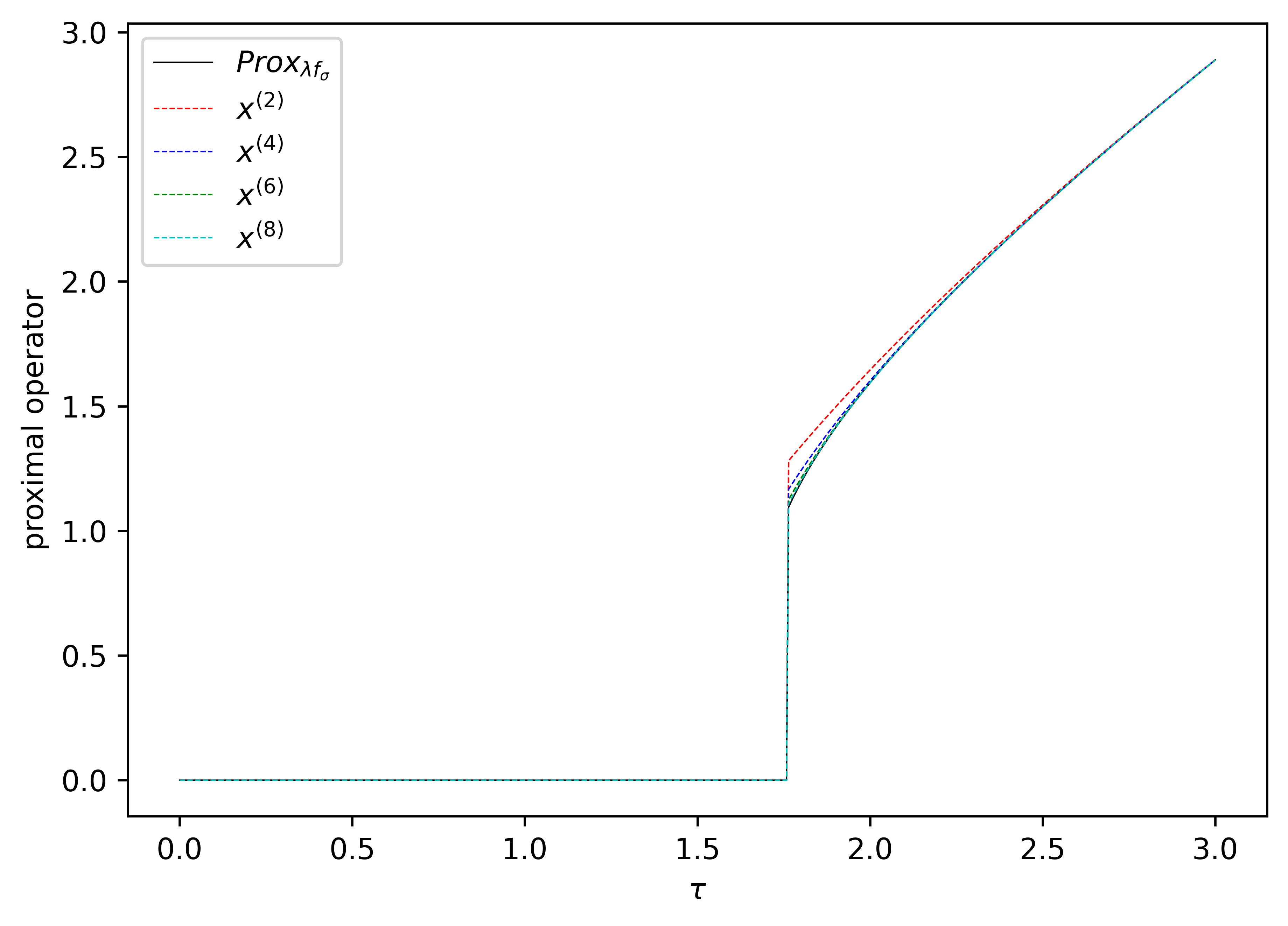}}
\subfloat[]{\includegraphics[width=0.45\textwidth]{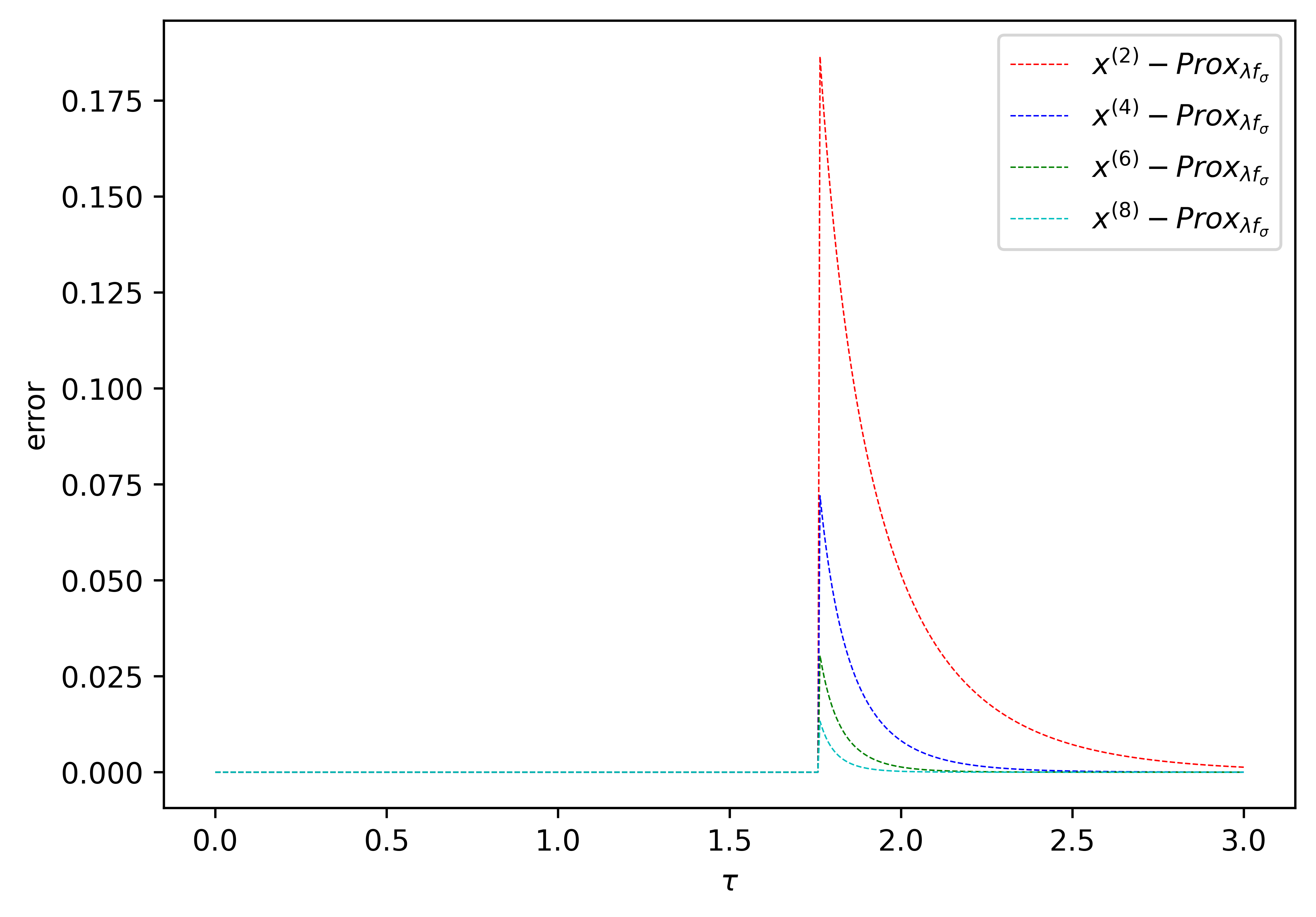}}
\caption{Illustration of Theorem \ref{Theorem4} with $\lambda=2$ and $\sigma=1$ for any $\tau\in(0,3]$. (a) The true proximal operator and IRL1 solution $x^{(k)}$ with $k=2,4,6,8$; (b) the error function.}
\label{PiE_prox_hard}
\end{figure}

Recall that $\sigma(1+\ln\frac{\lambda}{\sigma^2})\le \bar{\tau}_{\lambda,\sigma}\le \frac{\lambda}{\sigma}$ and $x_1(\tau)$ is strictly increasing for $\tau\ge\sigma(1+\ln\frac{\lambda}{\sigma^2})$. Observe that $x_1(\sigma(1+\ln\frac{\lambda}{\sigma^2}))=\sigma\ln\frac{\lambda}{\sigma^2}$.
It follows that for any $\tau>\bar{\tau}_{\lambda,\sigma}$, 
$$
\frac{\lambda}{\sigma^2}e^{-\frac{x_1(\tau)}{\sigma}}<\frac{\lambda}{\sigma^2}e^{-\frac{x_1(\bar{\tau}_{\lambda,\sigma})}{\sigma}}
\le \frac{\lambda}{\sigma^2}e^{-\frac{x_1(\sigma(1+\ln\frac{\lambda}{\sigma^2}))}{\sigma}}
=\frac{\lambda}{\sigma^2}e^{-\frac{\sigma\ln\frac{\lambda}{\sigma^2}}{\sigma}}
=1.
$$
Given the initial value $x^{(0)}$ as in Theorem \ref{Theorem4}. Let $\lambda=2$ and $\sigma=1$.
Fix $k\in\{2,4,6,8\}$, $x^{(k)}$, $\prox_{\lambda f_{\sigma}}$, and the corresponding error function $x^{(k)}-\prox_{\lambda f_{\sigma}}$ for any $\tau>0$ are given in Fig. \ref{PiE_prox_hard}.

\section{Proof of Theorems \ref{MainTh1} and \ref{MainTh2}}\label{Section4}

To start with, we present several technical lemmas describing the convergence of Algorithm \ref{Algorithm}. Then, the limit of the sequence generated by this algorithm is compared to $\prox_{\lambda f_{\sigma}}(\tau)$ directly with $\tau>0$.

\begin{lemma}\label{Exik0Lem} 
 Given $\tau\!>\! 0$ and $x^{(0)}\!\geq \!0$.
Let sequence $\{x^{(k)}\}$ be generated by Algorithm \ref{Algorithm}.   Suppose that there exists $k_0\!\ge\!0$ such that $x^{(k_0)}=0$. Then
\begin{itemize}
    \item [{\bf (i)}] If $\tau\le \frac{\lambda}{\sigma}$,  $x^{(k)}=0$ for any $k\geq k_0$.
    
    \item [{\bf (ii)}] If $\tau>\frac{\lambda}{\sigma}$,
     $ x^{(k+1)}=\tau-\frac{\lambda}{\sigma}e^{-\frac{x^{(k)}}{\sigma}}>0$ for any $k\geq k_0$. Moreover, $x^{(k+1)}> x^{(k)}$ for any $k\ge k_0$. 
    \item [{\bf (iii)}] 
    If $\tau>\frac{\lambda}{\sigma}$, the sequence $\{x^{(k)}\}$ converges to $\sigma W_0(-\frac{\lambda}{\sigma^2}e^{-\frac{\tau}{\sigma}})+\tau$.
\end{itemize}
\end{lemma}

\begin{proof} From \eqref{reweightedl1},  it holds that
\begin{equation}\label{reweightedl1eq1}
x^{(k_0+1)}=\Big(\tau-\frac{\lambda}{\sigma}e^{-\frac{x^{(k_0)}}{\sigma}}\Big)_+
=\Big(\tau-\frac{\lambda}{\sigma}\Big)_+.
\end{equation}
Clearly, if $\tau\!\le\! \frac{\lambda}{\sigma}$, $x^{(k_0+1)}\!=\!0$ by \eqref{reweightedl1eq1}, yielding  $x^{(k)}\!=\!0$ for any $k\!\ge\! k_0+1$. 
If $\tau\!>\!\frac{\lambda}{\sigma}$,  $x^{(k_0+1)}\!=\!\tau-\frac{\lambda}{\sigma}\!>\!0$ by \eqref{reweightedl1eq1}.  With  $\tau\!-\!\frac{\lambda}{\sigma}e^{-\frac{x^{(k_0+1)}}{\sigma}}\!>\!\tau\!-\!\frac{\lambda}{\sigma}e^{0}\!>\!0$ and \eqref{reweightedl1}, one has  
\[
x^{(k_0+2)}=(\tau-\frac{\lambda}{\sigma}e^{-\frac{x^{(k_0+1)}}{\sigma}})_+>0,
\]
 yielding 
\(
x^{(k+1)}=\tau-\frac{\lambda}{\sigma}e^{-\frac{x^{(k)}}{\sigma}}>0
\)
for any $k\ge k_0$.
Moreover, notice that $x^{(k_0+2)}>x^{(k_0+1)}$. Together with \(
x^{(k+1)}=\tau-\frac{\lambda}{\sigma}e^{-\frac{x^{(k)}}{\sigma}}>0
\)
for any $k\ge k_0$ and the monotonic increase of the function $h(t):=\tau-\frac{\lambda}{\sigma}e^{-t}
$, we know that
$x^{(k+1)}> x^{(k)}$ for any $k\ge k_0$.
Hence, the statements (i) and (ii) hold.

By (ii)  and $\tau-\frac{\lambda}{\sigma}\!<\!x^{(k)}\!<\!\tau$  for each $k>k_0$, the sequence $\{x^{(k)}\}$ converges. Moreover,
 it converges to $\sigma W_0(-\frac{\lambda}{\sigma^2}e^{-\frac{\tau}{\sigma}})+\tau$ by Proposition \ref{xinftyProp} (ii).
\end{proof}

  \begin{lemma}\label{MonoLem} 
 Given $\tau\!>\! 0$ and $x^{(0)}\!\geq \!0$.
Let sequence $\{x^{(k)}\}$ be generated by Algorithm \ref{Algorithm}.   If  $x^{(k)}>0$ for any $k\in \mathbb{N}$, then
\begin{itemize}
\item[{\bf(i)}] $\{x^{(k)}\}$ is strictly increasing and convergent  if $x^{(1)}>x^{(0)}$.
\item[{\bf(ii)}] $\{x^{(k)}\}$ is strictly decreasing and convergent  if $x^{(1)}<x^{(0)}$. 
\item[{\bf(iii)}] $\{x^{(k)}\}$ is constant if $x^{(1)}=x^{(0)}$. 
\end{itemize}
\end{lemma}
\begin{proof}
 Since $x^{(k)}\!>\!0$ for each $k$ and  \eqref{reweightedl1}, 
 $ x^{(k+1)}\!=\!\tau-\frac{\lambda}{\sigma}e^{-\frac{x^{(k)}}{\sigma}}\!>\!0$ for any $k\in \mathbb{N}$.
  If $x^{(1)}\!>\!x^{(0)}$, together with the monotonic increase of the function $h(t):=\tau-\frac{\lambda}{\sigma}e^{-t}
$, we know that
$x^{(k+1)}\!>\! x^{(k)}$ for any $k$.  
Obviously, $0<x^{(k)}<\tau$ for each $k$.
Hence, $\{x^{(k)}\}$ is strictly increasing and converging, and the statement (i) holds. The rest proof is similar to (i).
\end{proof}

\begin{lemma}\label{x0geq0Lem} 
Given $\tau\in [\frac{\lambda}{\sigma},+\infty)$ and $x^{(0)}>0$.
Let the sequence $\{x^{(k)}\}$ be generated by Algorithm \ref{Algorithm}. Then $x^{(k)}>0$ for all $k\ge0$ and the sequence $\{x^{(k)}\}$ converges to $x_1(\tau)$ defined as in Lemma  \ref{ProximalTheorem1}.
\end{lemma}
\begin{proof} 
 Since $\tau\!\geq\!\frac{\lambda}{\sigma}$ and $x^{(0)}\!>\!0$, it holds that $\tau-\frac{\lambda}{\sigma}e^{-\frac{x^{(0)}}{\sigma}}\!>\!\tau-\frac{\lambda}{\sigma}\!\geq\!0$.
 With \eqref{reweightedl1},  we have
\(
x^{(1)}=\Big(\tau-\frac{\lambda}{\sigma}e^{-\frac{x^{(0)}}{\sigma}}\Big)_+
=\tau-\frac{\lambda}{\sigma}e^{-\frac{x^{(0)}}{\sigma}}>0,
\)
which  yields 
\[
x^{(k+1)}=\tau-\frac{\lambda}{\sigma}e^{-\frac{x^{(k)}}{\sigma}}>0, \text{ for any } k\in\mathbb{N}.
\]
By Lemma \ref{MonoLem}, it suffices to argue the sequence $\{x^{(k)}\}$ converges to $x_1(\tau)$. 
Notice that
$x^{(1)}-x^{(0)}\!=\!\tau-\frac{\lambda}{\sigma}e^{-\frac{x^{(0)}}{\sigma}}-x^{(0)}\!=\!\phi(x^{(0)})$,
where
$\phi$ be defined by \eqref{phiDef}. Obviously, $x_1(\tau)\!\geq\! 0$ by Lemma \ref{phiLem} (iv) and (v). We will proceed in three cases.

{\bf Case 1: $x^{(0)}\!=\!x_1(\tau)\!>\!0$.} Then $\phi(x^{(0)})\!=\!0$ by Lemma \ref{phiLem} (iv), namely, $x^{(1)}\!=\!x^{(0)}$.  Hence,
$\{x^{(k)}\}$ is constant
from Lemma \ref{MonoLem} (iii).
The desired result obviously holds.

{\bf Case 2: $0\!<\!x^{(0)}\!<\!x_1(\tau)$}. Now, $x_1(\tau)\!>\!0$.  
Then $\phi(x^{(0)})\!>\!0$ by  Lemma \ref{phiLem} (i) 
 and the fact $\phi(0)\!\geq \! 0$, which implies that  $x^{(1)}\!>\!x^{(0)}$. Hence, $\{x^{(k)}\}$ is strictly increasing and convergent from Lemma \ref{MonoLem} (i).  

{\bf Case 3: $x^{(0)}\!>\!x_1(\tau)$}. Then $\phi(x^{(0)})\!<\!0$ by Lemma \ref{phiLem} (i), namely, $x^{(1)}\!<\!x^{(0)}$. Hence, $\{x^{(k)}\}$ is strictly decreasing and convergent from Lemma \ref{MonoLem} (ii).

In summary, the sequence $\{x^{(k)}\}$ is convergent and its limit is denoted by
$\!x^{(\infty)}$.  Then
 $x^{(\infty)}\!\geq\! 0$ and  $\phi(x^{(\infty)})\!=\!0$. So,
 $x^{(\infty)}=x_1(\tau)$ by Lemma \ref{phiLem} (iv) and (v).
\end{proof}

  By Lemma \ref{Exik0Lem} (iii) and  Lemma \ref{x0geq0Lem}, we have the following conclusion.

\begin{corollary}\label{CorLambSigm}
  Given $\tau\!>\! \frac{\lambda}{\sigma}$  and $x^{(0)}\!\geq\! 0$, 
  the sequence $\{x^{(k)}\}$ generated by Algorithm \ref{Algorithm} converges to $x_1(\tau)$. 
\end{corollary}

 The following lemma proves that  $\{x^{(k)}\}$ always converges to $0$ for all $\tau\in (0,\sigma(1+\ln \frac{\lambda}{\sigma^2}))$ if 
$\sigma(1+\ln \frac{\lambda}{\sigma^2})>0$, namely, $\frac{\lambda}{\sigma^2}>\frac{1}{e}$.

\begin{lemma}\label{xinftyCase2} Suppose $\sigma(1+\ln \frac{\lambda}{\sigma^2})>0$. Given $\tau\in(0,\sigma(1+\ln \frac{\lambda}{\sigma^2}))$ and an initial value $x^{(0)}\ge0$. Let the sequence $\{x^{(k)}\}$ be generated by Algorithm \ref{Algorithm}.  Then $\{x^{(k)}\}$ converges to $0$. 
\end{lemma}
\begin{proof} 
 Firstly, we will argue that there exists $k_0\ge0$ such that $x^{(k_0)}=0$.
 If not, $x^{(k)}>0$ for all $k\ge0$.
 Then $x^{(1)}\!=\!\tau-\frac{\lambda}{\sigma}e^{-\frac{x^{(0)}}{\sigma}}$ from \eqref{reweightedl1} and $x^{(1)}-x^{(0)}=\phi(x^{(0)})$, where $\phi$ be defined by \eqref{phiDef}.
Again from Lemma \ref{phiLem} (i) and 
$\tau\in(0,\sigma(1+\ln \frac{\lambda}{\sigma^2}))$, it holds that $\phi(x)\leq \phi(\sigma\ln{\frac{\lambda}{\sigma^2}})<0$ for any $x\in \mathbb{R}$. Consequently, $\phi(x^{(0)})<0$, namely, 
 $x^{(1)}<x^{(0)}$.
 So, $\{x^{(k)}\}$ is decreasing and convergent by Lemma \ref{MonoLem} (ii).
 Now suppose that $\lim\limits_{k\to \infty}x^{(k)}\!=\!x^{(\infty)}$. Then $x^{(\infty)}\!\geq \!0$ and $\phi(x^{(\infty)})\!=\!0$, which 
 contradicts to $\phi(x^{(\infty)})<0$. 
Hence, there exists $k_0\ge0$ such that $x^{(k_0)}=0$, and then the sequence $\{x^{(k)}\}$ converges to $0$ by Lemma \ref{Exik0Lem}  (i) and the fact $\sigma(1+\ln{\frac{\lambda}{\sigma^2}})\!\leq \!\frac{\lambda}{\sigma}$. 
\end{proof}

The next two lemmas study the convergence of $\{x^{(k)}\}$ for $\tau\!\in\! [\sigma(1+\ln \frac{\lambda}{\sigma^2}),\frac{\lambda}{\sigma})$.

\begin{lemma}\label{ComplexLemCase1} Suppose $\lambda\!\leq\! \sigma^2$. 
Given $\tau\!\in\! [\sigma(1+\ln \frac{\lambda}{\sigma^2}),\frac{\lambda}{\sigma})$ and an initial value $x^{(0)}\!\ge\!0$.
Let the sequence $\{x^{(k)}\}$ be generated by Algorithm \ref{Algorithm}. Then $\{x^{(k)}\}$ converges to $0$. 
\end{lemma}
\begin{proof} 
Firstly, we will argue that there exists $k_0\ge0$ such that $x^{(k_0)}=0$.
 If not, $x^{(k)}>0$ for all $k\ge0$.
 Then $x^{(1)}\!=\!\tau-\frac{\lambda}{\sigma}e^{-\frac{x^{(0)}}{\sigma}}$ from \eqref{reweightedl1} and $x^{(1)}-x^{(0)}=\phi(x^{(0)})$, where $\phi$ is defined by \eqref{phiDef}.
 Since $\lambda\le \sigma^2$, $\sigma\ln{\frac{\lambda}{\sigma^2}}\leq 0$.
Again from Lemma \ref{phiLem} (i), it holds that $\phi(x)\leq \phi(0)=\tau-\frac{\lambda}{\sigma}<0$ for any $x\geq 0$. Consequently, $\phi(x^{(0)})<0$, namely, 
 $x^{(1)}<x^{(0)}$.
 So, $\{x^{(k)}\}$ is decreasing and convergent by Lemma \ref{MonoLem} (ii).
 Now suppose that $\lim\limits_{k\to \infty}x^{(k)}\!=\!x^{(\infty)}$. Then $x^{(\infty)}\geq 0$ and $\phi(x^{(\infty)})=0$, which 
 contradicts to $\phi(x^{(\infty)})<0$. 
Hence, there exists $k_0\ge0$ such that $x^{(k_0)}=0$, and then the sequence $\{x^{(k)}\}$ converges to $0$ by Lemma \ref{Exik0Lem}  (i).
\end{proof}

\begin{lemma}\label{ComplexLemCase2} Suppose  $\lambda\!>\! \sigma^2$. 
Given $\tau\!\in\! [\sigma(1+\ln \frac{\lambda}{\sigma^2}),\frac{\lambda}{\sigma})$ and an initial value $x^{(0)}\!\ge\!0$.
Let the sequence $\{x^{(k)}\}$ be generated by Algorithm \ref{Algorithm} and 
$x_1(\tau)$ and $x_2(\tau)$ are defined in Lemma \ref{phiLem}. Then, the following statements hold.
\begin{itemize}
    \item [{\bf (i) }] If $x^{(0)}\in (0,x_2(\tau))$, the sequence $\{x^{(k)}\}$ converges to $0$.
    
    \item [{\bf (ii) }] If $x^{(0)}=x_2(\tau)$, the sequence $\{x^{(k)}\}$ converges to $x_2(\tau)$.
    
    \item [{\bf (iii) }]
    If $x^{(0)}\in (x_2(\tau),+\infty)$, the sequence $\{x^{(k)}\}$ converges to $x_1(\tau)$.
    \end{itemize}
\end{lemma}
\begin{proof} 
Since $\lambda>\sigma^2$ and $\tau\!\in\! [\sigma(1+\ln \frac{\lambda}{\sigma^2}),\frac{\lambda}{\sigma})$, 
$\phi(\sigma\ln\frac{\lambda}{\sigma\tau})\!=\!-\sigma\ln\frac{\lambda}{\sigma\tau}\!<\!0$, where $\phi$ is defined by \eqref{phiDef}. 
By Lemma \ref{phiLem} (i) and (ii), it holds
\begin{equation}\label{CompTempeq1}
0<\sigma\ln\frac{\lambda}{\sigma\tau}<x_2(\tau)\leq\sigma \ln \frac{\lambda}{\sigma^2}\leq x_1(\tau),
\end{equation}
for any $\tau\in [\sigma(1+\ln \frac{\lambda}{\sigma^2}),\frac{\lambda}{\sigma})$.  

(i) The proof can be divided into two cases: $x^{(0)}\le \sigma\ln\frac{\lambda}{\sigma\tau}$ and $\sigma\ln\frac{\lambda}{\sigma\tau}<x^{(0)}<x_2(\tau)$. 
If $x^{(0)}\le \sigma\ln\frac{\lambda}{\sigma\tau}$, then 
$$
x^{(1)}=(\tau-\frac{\lambda}{\sigma}e^{-\frac{x^{(0)}}{\sigma}})_+\le (\tau-\frac{\lambda}{\sigma}e^{-\frac{\sigma\ln\frac{\lambda}{\sigma\tau}}{\sigma}})_+=(\tau-\tau)_+=0,
$$
and hence $\{x^{(k)}\}$ converges to $0$ by Lemma \ref{Exik0Lem}  (i). 

If $\sigma\ln\frac{\lambda}{\sigma\tau}<x^{(0)}<x_2(\tau)$, 
\(
\tau-\frac{\lambda}{\sigma}e^{-\frac{x^{(0)}}{\sigma}}> \tau-\frac{\lambda}{\sigma}e^{-\frac{\sigma\ln\frac{\lambda}{\sigma\tau}}{\sigma}}=0.
\)
Hence,  it follows that $x^{(1)}=\tau-\frac{\lambda}{\sigma}e^{-\frac{x^{(0)}}{\sigma}}$ from \eqref{reweightedl1}, and 
then $0<x^{(1)}\le x^{(0)}$ as $x^{(1)}- x^{(0)}=\phi(x^{(0)})<\phi(x_2(
\tau))=0$ by Lemma \ref{phiLem} (i)--(iii). 
If there exists $k_0\ge0$ such that $x^{(k_0)}=0$,  $\{x^{(k)}\}$ converges to $0$ by Lemma \ref{Exik0Lem}  (i).
Otherwise,  $x^{(k)}>0$ for all $k\ge0$. 
Notice that $x^{(1)}<x^{(0)}$. Thus, the sequence $\{x^{(k)}\}$ is decreasing and convergent by Lemma \ref{MonoLem} (ii). Moreover, its limit, denoted by $x^{(\infty)}$  satisfies $x^{(\infty)}=\tau-\frac{\lambda}{\sigma}e^{-\frac{x^{(\infty)}}{\sigma}}$, namely, $\phi(x^{(\infty)})=0$, and $x^{(\infty)}<x^{(0)}<x_2(\tau)$, 
which implies  $\phi(x^{(\infty)})<\phi(x_2(\tau))=0$. Contradiction.
In summary, the sequence $\{x^{(k)}\}$ converges to $0$.

(ii) If $x^{(0)}\!=\!x_2(\tau)$, then $x^{(0)}\!>\!\sigma\ln\frac{\lambda}{\sigma\tau}$ from \eqref{CompTempeq1}, and $x^{(1)}\!=\!x^{(0)}$ as $x^{(1)}- x^{(0)}\!=\!\phi(x^{(0)})\!=\!\phi(x_2(\tau))\!=\!0$.
In this scenario, the sequence $\{x^{(k)}\}$ is a constant sequence and its limit is $x_2(\tau)$.

(iii) Let $x^{(0)}\!\in\!(x_2(\tau),+\infty)$.
 We know that $x_2(\tau)\!\leq\! x_1(\tau)$ from Lemma \ref{phiLem} (ii) and (iii).
 $x^{(0)}>\sigma\ln\frac{\lambda}{\sigma\tau}$ by \eqref{CompTempeq1} and
 \(
\tau-\frac{\lambda}{\sigma}e^{-\frac{x^{(0)}}{\sigma}}> \tau-\frac{\lambda}{\sigma}e^{-\frac{\sigma\ln\frac{\lambda}{\sigma\tau}}{\sigma}}=0.
\)
Hence,  it follows that $x^{(1)}=\tau-\frac{\lambda}{\sigma}e^{-\frac{x^{(0)}}{\sigma}}$ from \eqref{reweightedl1}.
If $x^{(0)}\!<\!x_1(\tau)$,
$x^{(1)}\!>\! x^{(0)}$ since $x^{(1)}- x^{(0)}\!=\!\phi(x^{(0)})\!>\!\phi(x_2(
\tau))\!=\!0$ by Lemma \ref{phiLem} (i) and (ii), which  yields 
\[
x^{(k+1)}=\tau-\frac{\lambda}{\sigma}e^{-\frac{x^{(k)}}{\sigma}}>0, \text{ for any } k\in \mathbb{N}.
\] 
 Hence, $\{x^{(k)}\}$ is increasing and convergent by Lemma \ref{MonoLem} (i), and its limit satisfies $x^{(\infty)}=\tau-\frac{\lambda}{\sigma}e^{-\frac{x^{(\infty)}}{\sigma}}$ and must be $x_1(\tau)$. 
 If $x^{(0)}=x_1(\tau)$, $x^{(1)}\!=\! x^{(0)}$ since $x^{(1)}- x^{(0)}=\phi(x^{(0)})=\phi(x_1(
\tau))=0$. In this scenario, the sequence $\{x^{(k)}\}$ is a constant sequence and its limit is $x_1(\tau)$.
If $x^{(0)}>x_1(\tau)$, 
 $x^{(1)}<x^{(0)}$  as $\phi(x^{(0)})<\phi(x_1(\tau))=0$ with Lemma \ref{phiLem} (i) and (ii). We estimate
$$
x^{(1)}-x_1(\tau)=\tau-\frac{\lambda}{\sigma}e^{-\frac{x^{(0)}}{\sigma}}-x_1(\tau)>\tau-\frac{\lambda}{\sigma}e^{-\frac{x_1(\tau)}{\sigma}}-x_1(\tau)=\phi(x_1(\tau))=0.
$$
which implies $x^{(1)}>x_1(\tau)$ and then $x^{(k)}>x^{(k+1)}>x_1(\tau)$ for each $k$. Therefore, $\{x^{(k)}\}$ is decreasing and convergent. Its limit satisfies $x^{(\infty)}\ge x_1(\tau)$ and $x^{(\infty)}=\tau-\frac{\lambda}{\sigma}e^{-\frac{x^{(\infty)}}{\sigma}}$. Therefore, $x^{(\infty)}$ must be $x_1(\tau)$ by Lemma \ref{phiLem} (ii) and (iii).
\end{proof}

By Lemma \ref{ComplexLemCase2} (ii), (iii) and Lemma \ref{phiLem} (iii), we  can obtain the following claim.
\begin {corollary}\label{CorTau}
 When $\tau=\sigma(1+\ln \frac{\lambda}{\sigma^2})$ and $\lambda>\sigma^2$,  the sequence $\{x^{(k)}\}$ converges to $x_1(\tau)$ for any $x^{(0)}\in [x_1(\tau),+\infty)$
 with $x_1(\tau)=\sigma \ln \frac{\lambda}{\sigma^2}$.
\end{corollary}

Now, we are ready to prove Theorems \ref{MainTh1} and \ref{MainTh2}.

{\bf Proof of Theorem \ref{MainTh1} } 
We only argue when $\tau\!>\!0$. In the following,
we will divide the arguments into two cases.

{\bf Case 1: $\sigma(1+\ln \frac{\lambda}{\sigma^2})\!>\!0$}. 
 When  $\tau\!\in\! (0,\sigma(1+\ln \frac{\lambda}{\sigma^2}))$,  $\{x^{(k)}\}$ converges to $0$ by Lemma \ref{xinftyCase2}.
When
$\tau\in [\sigma(1+\ln \frac{\lambda}{\sigma^2}),\frac{\lambda}{\sigma})$, $\{x^{(k)}\}$ converges to $0$ by Lemma \ref{ComplexLemCase1}. 
When $\tau=\frac{\lambda}{\sigma}$,  $\{x^{(k)}\}$ converges to $x_1(\tau)=0$ by Lemma \ref{x0geq0Lem} and Lemma \ref{phiLem} if $x^{(0)}>0$,  and
$\{x^{(k)}\}$ converges to $0$ by Lemma \ref{Exik0Lem} (i) if $x^{(0)}=0$. In a short, 
for  any  $\tau\in (0,\frac{\lambda}{\sigma}]$, $\{x^{(k)}\}$ converges to $0$. 
When $\tau\in (\frac{\lambda}{\sigma},+\infty)$,  $\{x^{(k)}\}$ converges to $x_1(\tau)$ by Lemma \ref{x0geq0Lem} if $x^{(0)}>0$, and
$\{x^{(k)}\}$ converges to $x_1(\tau)$ by Lemma \ref{Exik0Lem} (iii) if $x^{(0)}=0$. Hence, for any  $\tau\in (\frac{\lambda}{\sigma},+\infty)$, $\{x^{(k)}\}$ converges to $x_1(\tau)$.

{\bf Case 2: $\sigma(1+\ln \frac{\lambda}{\sigma^2})\leq 0$}. In this case.  
$\tau\in (0, \frac{\lambda}{\sigma})\subseteq [\sigma(1+\ln \frac{\lambda}{\sigma^2}),\frac{\lambda}{\sigma})$, the sequence $\{x^{(k)}\}$ converges to $0$ by Lemma \ref{ComplexLemCase1}. When  $\tau\in [\frac{\lambda}{\sigma},+\infty)$, its proof is the same as the case 1.

Based on the above arguments, $\{x^{(k)}\}$ converges to the exact solution to $\prox_{\lambda f_{\sigma}}(\tau)$ by Lemma \ref{ProximalTheorem1}. The proof is hence complete.

\bigskip
{\bf Proof of Theorem \ref{MainTh2} } 
We only argue that $\tau\!>\!0$.
By Corollary \ref{CorLambSigm}, $\{x^{(k)}\}$ converges to $x_1(\tau)$ for any $\tau\in(\frac{\lambda}{\sigma},+\infty)$. 
By Lemma \ref{xinftyCase2}, $\{x^{(k)}\}$ converges to $0$ for any $\tau\in (0,\sigma(1+\ln \frac{\lambda}{\sigma^2}))$. Hence, the statement (i) holds with Lemma \ref{ProximalTheorem2}.
The rest of the proof will focus on $\tau\in [\sigma(1+\ln \frac{\lambda}{\sigma^2}),\frac{\lambda}{\sigma}]$. Now suppose that $\tau\in [\sigma(1+\ln \frac{\lambda}{\sigma^2}),\frac{\lambda}{\sigma}]$.

(ii) Let $x^{(0)}\!\ge\!\sigma\ln \frac{\lambda}{\sigma^2}$.
If $\tau\in (\sigma(1+\ln \frac{\lambda}{\sigma^2}),\frac{\lambda}{\sigma}]$, then 
$x_2(\tau)\!<\! \sigma\ln \frac{\lambda}{\sigma^2}\!<\! x_1(\tau)$ by Lemma \ref{phiLem}(ii) and (v). Hence, $x^{(0)}\!>\!x_2(\tau)$ from the assumption that $x^{(0)}\!\ge\!\sigma\ln \frac{\lambda}{\sigma^2}$.  By Lemma \ref{ComplexLemCase2} (iii) and Corollary \ref{CorLambSigm}, $\{x^{(k)}\}$ converges to $x_1(\tau)$. If $\tau\!=\!\sigma(1+\ln \frac{\lambda}{\sigma^2})$,  $x_1(\tau)=x_2(\tau)=\sigma\ln \frac{\lambda}{\sigma^2}$ from  Lemma \ref{phiLem}(iii), and then the desired result is obtained by Corollary \ref{CorTau}. Thus, 
with Lemma \ref{ProximalTheorem3}, the statement (ii) holds. 

(iii) Let $x_2(\bar{\tau}_{\lambda,\sigma})\!<\!x^{(0)}\!<\! \sigma\ln\frac{\lambda}{\sigma^2}$.
Since  $x_2(\tau)$ is strictly decreasing on $\tau\!\in\! [\sigma(1+\ln \frac{\lambda}{\sigma^2}),\frac{\lambda}{\sigma})$ by Lemma \ref{WProperty} and $\sigma\ln\frac{\lambda}{\sigma^2}=x_2\Big(\sigma\big(1+\ln \frac{\lambda}{\sigma^2}\big)\Big)$ by Lemma \ref{phiLem}(iii),
we can drive that
\(
\sigma\Big(1+\ln \frac{\lambda}{\sigma^2}\Big)<x_2^{-1}(x^{(0)})<\bar{\tau}_{\lambda,\sigma},
\)
and that $x^{(0)}\!<\!x_2(\tau)$ for each $\tau\!\in\! [\sigma(1+\ln \frac{\lambda}{\sigma^2}),x_2^{-1}(x^{(0)})]$ and  $x^{(0)}\!>\!x_2(\tau)$ for each $\tau\!\in\! [x_2^{-1}(x^{(0)}),\frac{\lambda}{\sigma})$. Together with Lemma \ref{ComplexLemCase2}, the limit of $\{x^{(k)}\}$, denoted by $x^{(\infty)}$, satisfies 
\begin{equation}\label{limittauinfty}
x^{(\infty)}=\left\{\begin{array}{ll}
0,& \mbox{ if }\tau\in (\sigma(1+\ln \frac{\lambda}{\sigma^2}),x_2^{-1}(x^{(0)})),\\
x^{(0)}=x_2(\tau), &\mbox{ if  }\tau=x_2^{-1}(x^{(0)}),\\
x_1(\tau), &\mbox{ if  }\tau\in (x_2^{-1}(x^{(0)}),\frac{\lambda}{\sigma}].
\end{array}\right.
\end{equation}
Compared \eqref{limittauinfty} with Lemma \ref{ProximalTheorem3} gives the desired conclusion.

(iv) Let $x^{(0)}\!=\!x_2(\overline{\tau}_{\lambda,\sigma})$.
When $\tau\!\in\! [\sigma(1+\ln \frac{\lambda}{\sigma^2}),\overline{\tau}_{\lambda,\sigma}]$, $x_2(\tau)\!>\!x_2(\overline{\tau}_{\lambda,\sigma})\!=\!x_0$ since  $x_2(\tau)$ is strictly decreasing on $\tau\!\in\! [\sigma(1+\ln \frac{\lambda}{\sigma^2}),\frac{\lambda}{\sigma})$ by Lemma \ref{WProperty}, and then $\{x^{(k)}\}$ converges to $0$ by
 Lemma \ref{ComplexLemCase2} (i).
 When $\tau\!\in\! (\overline{\tau}_{\lambda,\sigma},\frac{\lambda}{\sigma}]$, $x_2(\tau)\!<\!x_2(\overline{\tau}_{\lambda,\sigma})\!=\!x_0$, and then $\{x^{(k)}\}$ converges to $x_1(\tau)$ by
 Lemma \ref{ComplexLemCase2} (iii).
  When $\tau\!=\!\overline{\tau}_{\lambda,\sigma}$, $x_2(\tau)\!=\!x_0$ and  then $\{x^{(k)}\}$ converges to $x_2(\overline{\tau}_{\lambda,\sigma})$ by Lemma \ref{ComplexLemCase2} (ii).
 Hence, the desired result is obtained  by Lemma \ref{ProximalTheorem3}.

(v) Let $0\le x^{(0)}\!<\!x_2(\bar{\tau}_{\lambda,\sigma})$. The proof is similar to (iii). Since $0\!\le x^{(0)}\!<\!x_2(\bar{\tau}_{\lambda,\sigma})$, we have
\(
\sigma(1+\ln \frac{\lambda}{\sigma^2})\!\leq\! \bar{\tau}_{\lambda,\sigma}\!<\!x_2^{-1}(x^{(0)})\!\le\! \frac{\lambda}{\sigma}.
\)
Suppose that $x^{(0)}\!>\!0$. Since  $x_2(\tau)$ is strictly decreasing on $\tau\!\in\! [\sigma(1+\ln \frac{\lambda}{\sigma^2}),\frac{\lambda}{\sigma}]$ by Lemma \ref{WProperty},
it holds that $x^{(0)}\!<\!x_2(\tau)$ for any $\tau\!\in\! [\sigma(1+\ln \frac{\lambda}{\sigma^2}),x_2^{-1}(x^{(0)}))$; and  $x^{(0)}\!>\!x_2(\tau)$ for any $\tau\!\in \!(x_2^{-1}(x^{(0)}),\frac{\lambda}{\sigma}]$.
By Lemma \ref{ComplexLemCase2}, the limit of $\{x^{(k)}\}$ is given as in \eqref{limittauinfty}.
Compared \eqref{limittauinfty} with Lemma \ref{ProximalTheorem3}, $x^{(\infty)}$ does not belong to $\prox_{\lambda f_{\sigma}}(\tau)$ for $\tau\!\in\! (\bar{\tau}_{\lambda,\sigma},x_2^{-1}(x^{(0)})]$. The rest is also true for $x^{(0)}\!=\!0$ by Lemma \ref{Exik0Lem} (i) and the facts that $x_2(\tau)\!=\!0$ if and only if $\tau\!=\!\frac{\lambda}{\sigma}$ from the proof of case (i) in Lemma \ref{phiLem} and 
$x_2(\tau)$ is strictly decreasing on $\tau\!\in\! [\sigma(1+\ln \frac{\lambda}{\sigma^2}),\frac{\lambda}{\sigma}]$.

\section{Conclusions}\label{Conclusions}

The relation between the IRL1 solution and the true proximal operator of PiE \eqref{PiE} has been clarified in Theorems \ref{MainTh1} and \ref{MainTh2}, which can be explicitly dependent upon $\sigma$, the initial value $x^{(0)}$, and the regularization parameter $\lambda$. Furthermore, to remedy the gap, the initial value was adaptively selected  as in Theorems \ref{Theorem3} and \ref{Theorem4} to guarantee that the IRL1  solution
belongs to the proximal operator of PiE. The results justify the usage of IRL1 for PiE whenever an initial value is appropriately given. Finally, our arguments can be applied to other sparse-promoting penalties, especially those whose proximal operator can not be explicitly derived.

\bibliographystyle{siam}
\bibliography{ProximalPiE}

\begin{thebibliography}{10}

\bibitem{Beck2017}
{\sc A.~Beck}, {\em First-order Methods in Optimization}, SIAM Publisher, Philadelphia, 2017.

\bibitem{Blumensath2008}
{\sc T.~Blumensath and M.~E. Davies}, {\em Iterative thresholding for sparse approximations}, Journal of Fourier Analysis and Applications, 14 (2008), pp.~629--654.

\bibitem{Bradley1998}
{\sc P.~S. Bradley and O.~L. Mangasarian}, {\em Feature selection via concave minimization and support vector machines}, in Proceedings of the 15th International Conference on Machine Learning, vol.~98, 1998, pp.~82--90.

\bibitem{Bradley1998b}
{\sc P.~S. Bradley, O.~L. Mangasarian, and W.~N. Street}, {\em Feature selection via mathematical programming}, INFORMS Journal on Computing, 10 (1998), pp.~209--217.

\bibitem{Candes2008}
{\sc E.~J. Candes, M.~B. Wakin, and S.~P. Boyd}, {\em Enhancing sparsity by reweighted $\ell_1$ minimization}, Journal of Fourier Analysis and Applications, 14 (2008), pp.~877--905.

\bibitem{Chen2014}
{\sc L.~Chen and Y.~Gu}, {\em The convergence guarantees of a non-convex approach for sparse recovery}, IEEE Transactions on Signal Processing, 62 (2014), pp.~3754--3767.

\bibitem{Fan2001}
{\sc J.~Fan and R.~Li}, {\em Variable selection via nonconcave penalized likelihood and its oracle properties}, Journal of the American Statistical Association, 96 (2001), pp.~1348--1360.

\bibitem{Fan2020}
{\sc J.~Fan, R.~Li, C.-H. Zhang, and H.~Zou}, {\em Statistical Foundations of Data Science}, Chapman and Hall/CRC, 2020.

\bibitem{Foucart2009}
{\sc S.~Foucart and M.-J. Lai}, {\em Sparsest solutions of underdetermined linear systems via $\ell_q$-minimization for $0<q\le1$}, Applied and Computational Harmonic Analysis, 26 (2009), pp.~395--407.

\bibitem{Fung2002}
{\sc G.~M. Fung, O.~L. Mangasarian, and A.~J. Smola}, {\em Minimal kernel classifiers}, Journal of Machine Learning Research, 3 (2002), pp.~303--321.

\bibitem{Gao2011}
{\sc C.~Gao, N.~Wang, Q.~Yu, and Z.~Zhang}, {\em A feasible nonconvex relaxation approach to feature selection}, in Proceedings of the AAAI Conference on Artificial Intelligence, vol.~25, 2011, pp.~356--361.

\bibitem{Guo2021}
{\sc W.~Guo, Y.~Lou, J.~Qin, and M.~Yan}, {\em A novel regularization based on the error function for sparse recovery}, Journal of Scientific Computing, 87 (2021), p.~No. 31.

\bibitem{Jiang2015}
{\sc W.~Jiang, F.~Nie, and H.~Huang}, {\em Robust dictionary learning with capped $\ell_1$-norm}, in Twenty-Fourth International Joint Conference on Artificial Intelligence, 2015.

\bibitem{Le2015}
{\sc H.~A. Le~Thi, T.~P. Dinh, H.~M. Le, and X.~T. Vo}, {\em {DC} approximation approaches for sparse optimization}, European Journal of Operational Research, 244 (2015), pp.~26--46.

\bibitem{LiuLin}
{\sc Y.~Liu and R.~Lin}, {\em A bisection method for computing the proximal operator of the $\ell_p$-norm with $0<p<1$}.
\newblock Journal of Computational and Applied Mathematics.

\bibitem{LiuZhouLin}
{\sc Y.~Liu, Y.~Zhou, and R.~Lin}, {\em The proximal operator of the piece-wise exponential function and its application in compressed sensing}.
\newblock arXiv:2306.13425.

\bibitem{Lou2018}
{\sc Y.~Lou and M.~Yan}, {\em Fast {L}1--{L}2 minimization via a proximal operator}, Journal of Scientific Computing, 74 (2018), pp.~767--785.

\bibitem{Lucidi2010}
{\sc S.~Lucidi and F.~Rinaldi}, {\em Exact penalty functions for nonlinear integer programming problems}, Journal of Optimization Theory and Applications, 145 (2010), pp.~479--488.

\bibitem{Malek2016}
{\sc M.~Malek-Mohammadi, A.~Koochakzadeh, M.~Babaie-Zadeh, M.~Jansson, and C.~R. Rojas}, {\em Successive concave sparsity approximation for compressed sensing}, IEEE Transactions on Signal Processing, 64 (2016), pp.~5657--5671.

\bibitem{Mangasarian1996}
{\sc O.~Mangasarian}, {\em Machine learning via polyhedral concave minimization}, in Applied Mathematics and Parallel Computing, Springer, 1996, pp.~175--188.

\bibitem{Mezo2022}
{\sc I.~Mezo}, {\em The Lambert W Function: Its Generalizations and Applications}, CRC Press, 2022.

\bibitem{Ochs2015}
{\sc P.~Ochs, A.~Dosovitskiy, T.~Brox, and T.~Pock}, {\em On iteratively reweighted algorithms for nonsmooth nonconvex optimization in computer vision}, SIAM Journal on Imaging Sciences, 8 (2015), pp.~331--372.

\bibitem{Prater2022}
{\sc A.~Prater-Bennette, L.~Shen, and E.~E. Tripp}, {\em The proximity operator of the log-sum penalty}, Journal of Scientific Computing, 93 (2022), p.~No. 67.

\bibitem{Prater2023}
{\sc A.~Prater-Bennette, L.~Shen, and E.~E. Tripp}, {\em A constructive approach for computing the proximity operator of the p-th power of the $\ell_1$ norm}, Applied and Computational Harmonic Analysis, 67 (2023), p.~101572.

\bibitem{Rinaldi2009}
{\sc F.~Rinaldi}, {\em New results on the equivalence between zero-one programming and continuous concave programming}, Optimization Letters, 3 (2009), pp.~377--386.

\bibitem{Rockafellar2009}
{\sc R.~T. Rockafellar and R.~J.-B. Wets}, {\em Variational analysis}, vol.~317, Springer Science \& Business Media, 2009.

\bibitem{Shen2016}
{\sc L.~Shen, Y.~Xu, and X.~Zeng}, {\em Wavelet inpainting with the $\ell_0$ sparse regularization}, Applied and Computational Harmonic Analysis, 41 (2016), pp.~26--53.

\bibitem{Tao2022}
{\sc M.~Tao}, {\em Minimization of ${L}_1$ over ${L}_2$ for sparse signal recovery with convergence guarantee}, SIAM Journal on Scientific Computing, 44 (2022), pp.~A770--A797.

\bibitem{Trzasko2008}
{\sc J.~Trzasko and A.~Manduca}, {\em Highly undersampled magnetic resonance image reconstruction via homotopic $\ell_0$-minimization}, IEEE Transactions on Medical imaging, 28 (2008), pp.~106--121.

\bibitem{Wang2023}
{\sc H.~Wang, H.~Zeng, and J.~Wang}, {\em Convergence rate analysis of proximal iteratively reweighted $\ell_1$ methods for $\ell_p$ regularization problems}, Optimization Letters, 17 (2023), pp.~413--435.

\bibitem{WangZeng2021}
{\sc H.~Wang, H.~Zeng, J.~Wang, and Q.~Wu}, {\em Relating $\ell_p$ regularization and reweighted $\ell_1$ regularization}, Optimization Letters, 15 (2021), pp.~2639--2660.

\bibitem{Wang2021}
{\sc H.~Wang, F.~Zhang, Y.~Shi, and Y.~Hu}, {\em Nonconvex and nonsmooth sparse optimization via adaptively iterative reweighted methods}, Journal of Global Optimization, 81 (2021), pp.~717--748.

\bibitem{Wright2022}
{\sc J.~Wright and Y.~Ma}, {\em High-dimensional Data Analysis with Low-dimensional Models: Principles, Computation, and Applications}, Cambridge University Press, 2022.

\bibitem{Yan2022}
{\sc J.~Yan, X.~Meng, F.~Cao, and H.~Ye}, {\em A universal rank approximation method for matrix completion}, International Journal of Wavelets, Multiresolution and Information Processing, 20 (2022), p.~2250016.

\bibitem{Yin2014}
{\sc P.~Yin, E.~Esser, and J.~Xin}, {\em Ratio and difference of $\ell_1$ and $\ell_2$ norms and sparse representation with coherent dictionaries}, Communications in Information and Systems, 14 (2014), pp.~87--109.

\bibitem{Yin2015}
{\sc P.~Yin, Y.~Lou, Q.~He, and J.~Xin}, {\em Minimization of $\ell_{1-2}$ for compressed sensing}, SIAM Journal on Scientific Computing, 37 (2015), pp.~A536--A563.

\bibitem{Zhang2010}
{\sc C.-H. Zhang}, {\em Nearly unbiased variable selection under minimax concave penalty}, Annals of Statistics, 38 (2010), pp.~894--942.

\bibitem{Zhang2017}
{\sc S.~Zhang and J.~Xin}, {\em Minimization of transformed $l_1$ penalty: Closed form representation and iterative thresholding algorithms}, Communications in Mathematical Sciences, 15 (2017), pp.~511--537.

\bibitem{Zhang2018}
\leavevmode\vrule height 2pt depth -1.6pt width 23pt, {\em Minimization of transformed ${L}_1$ penalty: theory, difference of convex function algorithm, and robust application in compressed sensing}, Mathematical Programming, 169 (2018), pp.~307--336.

\bibitem{Zhang2010b}
{\sc T.~Zhang}, {\em Analysis of multi-stage convex relaxation for sparse regularization}, Journal of Machine Learning Research, 11 (2010), pp.~1081--1107.

\bibitem{Zhou2022}
{\sc Z.~Zhou}, {\em A unified framework for constructing nonconvex regularizations}, IEEE Signal Processing Letters, 29 (2022), pp.~479--483.

\bibitem{Zhou2023}
\leavevmode\vrule height 2pt depth -1.6pt width 23pt, {\em Sparse recovery based on the generalized error function}, Journal of Computational Mathematics,doi:10.4208/jcm.2204-m2021-0288,  (2023).

\bibitem{Zou2008}
{\sc H.~Zou and R.~Li}, {\em One-step sparse estimates in nonconcave penalized likelihood models}, Annals of Statistics, 36 (2008), pp.~1509--1533.

\end{thebibliography}
\end{document}